\documentclass[a4paper,11pt,twoside]{article}
\usepackage{pst-fill,pst-grad}
\usepackage{textcomp}
\usepackage[english]{babel}
\usepackage{graphicx}
\usepackage{amsmath}
\usepackage{float}
\usepackage[matrix,arrow,curve]{xy}
\usepackage{pstricks}
\usepackage{amsmath,amsfonts,verbatim,afterpage,theorem,euscript,mathrsfs,amssymb}
\usepackage{amsfonts}
\usepackage{amssymb}
\usepackage{array}
\usepackage{dsfont}
\usepackage[colorlinks,citecolor=blue,hypertexnames=false,urlcolor=blue, linkcolor=blue, breaklinks=true]{hyperref}
\newtheorem{Theoreme}{Theorem}

\newtheorem{Proposition}{Proposition}[section]
\newtheorem{Lemme}{Lemma}[section]

\newtheorem{Remarque}{Remark}[section]
\newcommand{\mysection}{\setcounter{equation}{0} \section}

\def\vn{\vec{\nabla}}
\def\R{\mathbb{R}^n}

\title{\bf The role of the dimension in uniqueness results for the stationary quasi-geostrophic system} 
\usepackage{authblk}
\author{Diego Chamorro\footnote{\emph{diego.chamorro@univ-evry.fr} (corresponding author) }}
\affil{\footnotesize LaMME, Univ. Evry, CNRS, Universit\'e Paris-Saclay, 91025, Evry, France.}
\author{Manuel Fernando  Cortez\footnote{\emph{manuel.cortez@epn.edu.ec}}}
\affil{\footnotesize Departamento de matem\'aticas, Escuela Polit\'ecnica Nacional, Quito, Ecuador.}
\begin{document}
\maketitle
\begin{scriptsize}
\abstract{In this paper, we study a Liouville-type theorem for the stationary fractional quasi-geostrophic equation in various dimensions. Indeed, our analysis focuses on dimensions $n=2, 3, 4$ and we explore the uniqueness of weak solutions for this fractional system. We demonstrate here that, under some specific Lebesgue integrability information, the only admissible solution to the stationary fractional quasi-geostrophic system is the trivial one and this result provides a comprehensive understanding of how the dimension in connection to the fractional power of the Laplacian influences the uniqueness properties of weak solutions.}\\[3mm]
\textbf{Keywords: quasi-geostrophic equation; Liouville type theorems; uniqueness.}\\
\textbf{Mathematics Subject Classification: 76D03; 35A02; 35B65.} 
\end{scriptsize}
\mysection{Introduction and presentation of the results}

The study of fractional equations has gained significant attention in recent years due to their ability to model complex physical phenomena involving turbulence and fluid dynamics in various contexts, see \emph{e.g.} \cite{Constantin}, \cite{Pedlosky}, \cite{Smith}, \cite{Su}. One of the fundamental questions in the theory of partial differential equations (PDEs) is the existence, uniqueness and regularity of solutions to these equations and in this regard, Liouville-type theorems may play a crucial role by establishing conditions under which a solution to some equation is trivial (\emph{i.e.} identically zero).\\ 

In particular, this paper focuses on the stationary fractional quasi-geostrophic equation which is used to describe phenomena in atmospheric and oceanic dynamics (for more details see \cite{Constantin1} and \cite{Resnick}) and we will establish some uniqueness results for the following system:
\begin{equation}\label{SQGStationnaire}
(-\Delta)^{\frac \alpha 2} \theta +\mathbb{A}_{[\theta]}\cdot \vn \theta-f=0,\quad\mbox{with } 0<\alpha<2, \quad x\in \R,
\end{equation}
where $\theta: \R\longrightarrow \mathbb{R}$ is the variable,  the fractional operator $(-\Delta)^{\frac \alpha 2}$ is defined at the Fourier level by the symbol $|\xi|^\alpha$ and $f: \R\longrightarrow \mathbb{R}$ is a given external force. The vector $\mathbb{A}_{[\theta]}$ depends on $\theta$ in the following manner:
$$\mathbb{A}_{[\theta]}=[A_1(\theta), A_2(\theta),\cdots, A_n(\theta)],$$ 
where the operators $A_j$ for $1\leq j\leq n$ are assumed to be singular integral operators that satisty the boundedness property
\begin{equation}\label{OperateurAContinuite}
\|\mathbb{A}_{[\theta]}\|_{L^p}\leq C_{\mathbb{A}}\|\theta\|_{L^p},
\end{equation}
for all $1<p<+\infty$ (as long as $\theta \in L^p(\R)$). Moreover, the vector $\mathbb{A}_{[\theta]}$ is assumed to be divergence-free, \emph{i.e.} we have 
\begin{equation}\label{OperateurADivergence}
div(\mathbb{A}_{[\theta]})=0.
\end{equation}
Although the evolution variants of the previous system was considerably studied in \cite{Cordoba}, \cite{Dong} and \cite{Michael} (see also the references therein), to the best of our knowledge, the study of the stationary system (\ref{SQGStationnaire}) does not seem to have been done (see however \cite{Bahrami} for the inviscid 2D problem).\\

Thus, the main purpose of this article is to establish the uniqueness of the trivial solution for the system (\ref{SQGStationnaire}) under some hypothesis stated in terms of Lebesgue spaces. This research program is known in the literature as \emph{Liouville-type theorems} and has been applied mainly for the Navier-Stokes equations, see \cite{ChaeYoneda}, \cite{ChaeWolf}, \cite{Galdi},  \cite{OJ1}, \cite{Koch}, \cite{Kozono}, \cite{Ser2015} and the references therein for more details. See also \cite{wx18} for the fractional Navier-Stokes problem.\\

To start, in our first result we construct weak solutions for the system (\ref{SQGStationnaire}) in the space $\dot{H}^{\frac{\alpha}{2}}(\R)$:
\begin{Theoreme}[Existence]\label{Theo_ExistenceStationnaires}
Consider $n\geq 2$. Fix $0< \alpha< 2$ and consider $f\in \dot{H}^{-\frac{\alpha}{2}}(\R)$ an external force. There exists  $\theta\in\dot{H}^{\frac \alpha 2}(\R)$ that satisfies in the weak sense the equation (\ref{SQGStationnaire}).
\end{Theoreme}
\noindent The proof of this theorem is rather standard and it relies in an application of the Leray-Schauder fixed point argument (see \cite[Chapter 16]{PGLR2}). Note in particular that the resolution  space $\dot{H}^{\frac \alpha 2}(\R)$ for the system (\ref{SQGStationnaire}) is naturally determined by the structure of the equation when looking for energy estimates. Indeed, if we (formally) multiply the equation (\ref{SQGStationnaire}) by $\theta$ and we integrate, we obtain
$$\int_{\R}[(-\Delta)^{\frac \alpha 2} \theta +\mathbb{A}_{[\theta]}\cdot \vn \theta-f]\theta dx=0,$$
using the divergence free property of the vector field $\mathbb{A}$, we easily obtain the expression
$$\int_{\R}|(-\Delta)^{\frac{\alpha}{4}}\theta|^2dx=\int_{\R}f\theta dx,$$
from which, by the $\dot{H}^{-\frac{\alpha}{2}}-\dot{H}^{\frac{\alpha}{2}}$ duality and by the Young inequalities for the product, we deduce the control 
$$\|\theta\|_{\dot{H}^{\frac{\alpha}{2}}(\R)}\leq C\|f\|_{\dot{H}^{-\frac{\alpha}{2}}(\R)}.$$
Although these computations are completely formal (we will provide a rigorous proof in the Appendix \ref{Secc_ProofTheo1}), they give a good idea why the space $\dot{H}^{\frac{\alpha}{2}}(\R)$ natural. Remark however that the existence results are obtained via an approximation procedure and the Leray-Schauder fixed point argument does not provide any information about the uniqueness of the solutions.\\

In the following result, which is the main theorem of this article, we will address the problem of uniqueness of the weak solutions of the system (\ref{SQGStationnaire}) when $f=0$. As we shall see, the behavior of this problem will vary with the dimension $n=2,3,4$ and with the values of the fractional power of the Laplacian $0<\alpha<2$. Let us mention that the key ingredient of our analysis is related to a suitable decay at infinity of the solutions $\theta$. Indeed, due to the Sobolev embeddings, we have $\theta \in \dot{H}^{\frac{\alpha}{2}}(\R)\subset L^{\frac{2n}{n-\alpha}}(\R)$ so we dispose of a ``natural'' decay at infinity of the solution $\theta$, but this information doesn't seem enough to deduce that the unique solution of the equation  (\ref{SQGStationnaire}) when $f=0$ is the trivial solution $\theta \equiv 0$.  However, under some suitable additional information, stated below in terms of Lebesgue spaces, then we can obtain the uniqueness of the trivial solution:

\begin{Theoreme}[Uniqueness]\label{Theo_LiouvilleSQG} In the setting of the Theorem \ref{Theo_ExistenceStationnaires} with $n=2,3,4$, consider a solution $\theta \in \dot{H}^{\frac{\alpha}{2}}(\R)$ of the equation (\ref{SQGStationnaire}) and assume that $f=0$. Assume that is $\theta$ smooth and we have $\theta\in L^{\frac{2n-\epsilon}{n-\alpha}}(\R)$ with $0<\epsilon\ll 1$, moreover:
\begin{itemize}
\item[1)]  If $0<\alpha<\frac{n}{3}<2$, for $n=2,3,4$ we will ask the additional condition $\theta \in L^{\frac{2n+\nu}{n-\alpha}}(\R)$ for some real parameter $\nu$ such that  $(n-3\alpha)<\nu <1+(n-3\alpha)$,
\item[2)] If $\frac{n+2}{3}\leq \alpha<2$ for $n=2,3$, we will assume $\theta \in L^{\frac{3n}{n-1}}(\R)$,
\end{itemize}
then, under these assumptions, the unique solution of the system (\ref{SQGStationnaire}) is the trivial solution $\theta \equiv 0$.
\end{Theoreme}
\noindent Some remarks are in order there. First note that since we are considering $\dot{H}^{\frac{\alpha}{2}}(\R)$ solutions we have at our disposal (via the Sobolev embeddings) the critical Lebesgue information $L^{\frac{2n}{n-\alpha}}(\R)$. However, and regardless to the dimension, this information seems not enough when dealing with the fractional Laplacian: indeed, at some point (see the formula (\ref{ApplicationKatoPonce}) below) we need to use the Leibniz fractional rule to study terms of the form $(-\Delta)^{\frac{\alpha}{2}}(\theta \varphi_R)$ (where $\theta$ is the solution of the equation (\ref{SQGStationnaire}) and $\varphi_R$ is a cut off function) and this requires some additional information which is given by the hypothesis $\theta\in L^{\frac{2n-\epsilon}{n-\alpha}}(\R)$  -which is close to the critical Lebesgue space since the positive parameter $\epsilon>0$ can be made very small but not equal to $0$. Note next that outside the two cases stated in the theorem above, it is not necessary to impose hypothesis in order to obtain the uniqueness of the trivial solution. Remark also that the conditions $1)$ and $2)$ asked in the theorem above are essentially technical and they are related to the algebraic relationships  between the indexes of the Sobolev embeddings. Note that we only considered here the dimensions $n=2,3,4$ but the same ideas and techniques can be extended without problem to higher dimensions. Let us point out finally that the uniqueness of the trivial solution $\theta\equiv 0$ with the sole information $\theta\in \dot{H}^{\frac{\alpha}{2}}(\R)\subset L^{\frac{2n}{n-\alpha}}(\R)$ remains, to the best of our knowledge, a difficult open problem. \\

In the Theorem \ref{Theo_LiouvilleSQG} above we assumed that the solutions are smooth. It is worth to mention here that the problem of the regularity of the weak solutions for the system (\ref{SQGStationnaire}) is naturally related to the power $0<\alpha<2$ of the fractional Laplacian. Indeed, if this parameter $\alpha$ is big enough then the regularity of the solutions is not difficult to obtain (see the Appendix \ref{Secc_Regularity} below), however if $\alpha$ is small (in particular when $0<\alpha\leq 1$), then the regularity of the weak solutions obtained via the Theorem \ref{Theo_ExistenceStationnaires} is a completely open problem which is not studied here. \\

The plan of the article is the following. In Section \ref{Secc_Notation} we recall some of the tools that will be used in this article and in Section \ref{Secc_Uniqueness} we prove the main result (Theorem \ref{Theo_LiouvilleSQG}). In the Appendix \ref{Secc_ProofTheo1} we prove the Theorem \ref{Theo_ExistenceStationnaires} and in the Appendix \ref{Secc_Regularity} we establish some regularity results for the system (\ref{SQGStationnaire}). 
\mysection{Preliminaries}\label{Secc_Notation}
For $1<p<+\infty$ and for $s>0$ we define the homogeneous Sobolev spaces $\dot{W}^{s,p}(\R)$ by the condition 
$$\|f\|_{\dot{W}^{s,p}(\R)}=\|(-\Delta)^{\frac{s}{2}}f\|_{L^p(\R)}<+\infty.$$
In the special case when $p=2$ we simply write $\dot{W}^{s,2}(\R)=\dot{H}^{s}(\R)$. The non-homogeneous Sobolev spaces $W^{s,p}(\R)$ are defined by the condition 
$$\|f\|_{W^{s,p}(\R)}=\|f\|_{L^p(\R)}+\|(-\Delta)^{\frac{s}{2}}f\|_{L^p(\R)}<+\infty,$$
from which we easily deduce the embedding $W^{s,p}(\R)\subset\dot{W}^{s,p}(\R)$. Note also that, if $s_1>s_0>0$ then we have the space inclusion $W^{s_1,p}(\R)\subset W^{s_0,p}(\R)$. As the Sobolev spaces will constitute our main framework, we recall in the following lemmas some classical and useful results. 
\begin{Lemme}[Sobolev embeddings]\label{Lemm_SoboIneq}
\begin{itemize}
\item[]
\item[1)]For $0<s<\frac{n}{p}$ and $1<p,q<+\infty$, if we have the relationship $-\frac{n}{q}=s-\frac{n}{p}$, then we have the classical Sobolev inequality
$$\|f\|_{L^q(\R)}\leq C\|f\|_{\dot{W}^{s,p}(\R)}.$$
In particular we will use 
$$\|f\|_{L^{\frac{2n}{n-\alpha}}(\R)}\leq C\|f\|_{\dot{H}^{\frac{\alpha}{2}}(\R)}.$$
\item[2)] If $s=\frac{n}{p}$, then we have the space embedding $W^{s,p}(\R)\subset L^{q}(\R)$ for all $p\leq q<+\infty$, \emph{i.e.} we have
$$\|f\|_{L^q(\R)}\leq C\|f\|_{\dot{W}^{s,p}(\R)}.$$
\item[3)]If $0<s_0<s_1$ and $1<p_0, p_1<+\infty$ are such that $s_0-\frac{n}{p_0}=s_1-\frac{n}{p_1}$, then we have the following Sobolev space inclusion:
$$\dot{W}^{s_1,p_1}(\R)\subset\dot{W}^{s_0,p_0}(\R).$$
\end{itemize}
\end{Lemme}
\begin{Lemme}[Fractional Leibniz rule]\label{FracLeibniz}
\begin{itemize}
\item[]
\item[1)] Consider $f,g$ two smooth functions. Then we have the estimate
$$\|(-\Delta)^{\frac{s}{2}}(fg)\|_{L^p(\R)}\leq C\|(-\Delta)^{\frac{s}{2}}f\|_{L^{p_0}(\R)}\|g\|_{L^{p_1}(\R)}+C\|f\|_{L^{q_0}}\|(-\Delta)^{\frac{s}{2}}g\|_{L^{q_1}(\R)},$$
where $\frac1p=\frac{1}{p_0}+\frac{1}{p_1}=\frac{1}{q_0}+\frac{1}{q_1}$, with $0<s$, $1<p<+\infty$ and $1<p_0,p_1, q_0, q_1\leq +\infty$.
\item[2)] For $0<s, s_1, s_2<1$ with $s=s_1+s_2$ and $1<p, p_1, p_2<+\infty$ with $\frac{1}{p}=\frac{1}{p_1}+\frac{1}{p_2}$, we have
$$\|(-\Delta)^{\frac{s}{2}}(fg)-(-\Delta)^{\frac{s}{2}}(f)g-(-\Delta)^{\frac{s}{2}}(g)f\|_{L^p(\R)}\leq C\|(-\Delta)^{\frac{s_1}{2}}f\|_{L^{p_1}(\R)}\|(-\Delta)^{\frac{s_2}{2}}g\|_{L^{p_2}(\R)}.$$
\end{itemize}
\end{Lemme}
See \cite{Naibo} and \cite{GrafakosOh} for a proof of these estimates.\\

In the case of the $L^2$-based Sobolev spaces we also have the following useful estimate: 
\begin{Lemme}[Product rule in Sobolev spaces]\label{ProductRule}
For $0\leq s<+\infty$ and $0<\delta<\frac{n}{2}$,
$$\|fg\|_{\dot{H}^{s+\delta-\frac{n}{2}}(\R)}\leq C\left(\|f\|_{\dot{H}^{\delta}(\R)}\|g\|_{\dot{H}^{s}(\R)}+\|g\|_{\dot{H}^{\delta}(\R)}\|f\|_{\dot{H}^{s}(\R)}\right).$$
\end{Lemme}
See \cite[Lemma 7.3]{PGLR2} for a proof of this inequality.

\mysection{Uniqueness}\label{Secc_Uniqueness}
We study now the announced uniqueness results for the equation (\ref{SQGStationnaire}) and we will start with some steps that are independent from the dimension. Indeed, consider $\theta \in \dot{H}^{\frac{\alpha}{2}}(\R)$ a smooth solution of the equation (\ref{SQGStationnaire}) and for $R>1$ let us introduce the cut off function $\varphi_R$ defined by
\begin{equation}\label{Def_CutoffTest}
\varphi_R(x) = \varphi\left(\frac{x}{R}\right),
\end{equation}
where $\varphi \in \mathcal{C}_0^\infty(\mathbb{R}^n)$, with $\varphi(x) = 1$ for $|x| \leq \frac{1}{2}$ and $\varphi(x) = 0$ for $|x| \geq 1$. We gather in the next lemma some properties of the function $\varphi_R$ given above: 
\begin{Lemme}\label{LemmaN2}
For the function $\varphi_R$ defined in (\ref{Def_CutoffTest}) we have:
\begin{itemize}
\item[1)] the support of $\varphi_R$ satisfies $supp(\varphi_R)\subset B_R=B(0,R)$,
\item[2)] for all $1\leq p\leq +\infty$ we have $\|\varphi_R\|_{L^p(\R)} \leq C R^{\frac{n}{p}}$,
\item[3)] the support of $\nabla \varphi_R$ and $\Delta \varphi_R$ satisfies 
$$supp(\nabla \varphi_R), supp(\Delta \varphi_R)\subset \mathcal{C}_R = \{ x \in \mathbb{R}^n : \tfrac{R}{2} \leq |x| \leq R \},$$
\item[4)] for all $1\leq p\leq +\infty$ we have $\|\nabla \varphi_R\|_{L^p(\R)}=\|\nabla \varphi_R\|_{L^p(\mathcal{C}_R)} \leq C R^{-1 + \frac{n}{p}}$ and we also have $\|\Delta \varphi_R\|_{L^p(\R)} =\|\Delta \varphi_R\|_{L^p(\mathcal{C}_R)} \leq C R^{-2 + \frac{n}{p}}$,
\item[5)] for all $1<p<+\infty$ and for $\alpha>0$ we have $\|(-\Delta)^{\frac{\alpha}{2}} \varphi_R\|_{L^p(\R)} \leq C R^{-\alpha + \frac{n}{p}} \).
\end{itemize}
\end{Lemme}
The proof of these facts is straightforward and left to the reader.\\

\noindent Now we multiply the equation (\ref{SQGStationnaire}) by $\theta\varphi_R$ and we integrate to obtain (recall that we have assumed here that $f=0$ and that the solution $\theta$ is smooth):
$$\int_{\R}\left[(-\Delta)^{\frac \alpha 2} \theta +\mathbb{A}_{[\theta]}\cdot \vn \theta\right](\theta \varphi_R)dx=0.$$
Due to the divergence free property of the vector $\mathbb{A}_{[\theta]}$ we have the identity $\mathbb{A}_{[\theta]}\cdot \vn \theta=div(\mathbb{A}_{[\theta]} \theta)$ and we write
\begin{equation}\label{Identite1}
\int_{\R}\left((-\Delta)^{\frac \alpha 2} \theta\right)(\theta \varphi_R)dx+\int_{\R}div(\mathbb{A}_{[\theta]} \theta)(\theta \varphi_R)dx=0.
\end{equation} 
We study now each one of these terms separately.\\
\begin{itemize}
\item Indeed, for the first term in (\ref{Identite1}) we write, using the properties of the operator $(-\Delta)^{\frac{\alpha}{2}}$:
\begin{eqnarray*}
\int_{\R} \left((-\Delta)^{\frac{\alpha}{2}} \theta\right) \left(\theta\varphi_R\right)dx&=&\int_{\R} (-\Delta)^{\frac{\alpha}{4}}\theta\left[  (-\Delta)^{\frac{\alpha}{4}}\left(\theta\varphi_R\right)\right]dx\\
&=&\int_{\R} (-\Delta)^{\frac{\alpha}{4}} \theta \left[\left((-\Delta)^{\frac{\alpha}{4}}\theta\right) \varphi_R + (-\Delta)^{\frac{\alpha}{4}}( \theta \varphi_R)-\left((-\Delta)^{\frac{\alpha}{4}}\theta\right) \varphi_R\right]dx,
\end{eqnarray*}
and we have
\begin{eqnarray}
\int_{\R} \left((-\Delta)^{\frac{\alpha}{2}} \theta\right) \left(\theta\varphi_R\right)dx
&=&\int_{B_R} |(-\Delta)^{\frac{\alpha}{4}} \theta|^2\varphi_R dx\notag\\
&&+ \int_{\R} (-\Delta)^{\frac{\alpha}{4}} \theta\left[(-\Delta)^{\frac{\alpha}{4}}(\theta\varphi_R)-\left((-\Delta)^{\frac{\alpha}{4}}\theta\right) \varphi_R\right]dx,\label{Identite2}
\end{eqnarray}
where we used the fact that $supp(\varphi_R)\subset B_R$ in the second integral above. \\

\item For the second term of (\ref{Identite1}) we have by an integration by parts:
\begin{eqnarray*}
\int_{\R}div(\mathbb{A}_{[\theta]} \theta)(\theta \varphi_R)dx&=&-\int_{\R} \left(\mathbb{A}_{[\theta]}\cdot \vn \theta\right)\left(\theta\varphi_R\right)dx-\int_{\R} \theta(\mathbb{A}_{[\theta]}\theta)\cdot\vn\varphi_R dx\\
&=&-\int_{\R} div(\mathbb{A}_{[\theta]} \theta)(\theta \varphi_R)dx-\int_{\R} \theta(\mathbb{A}_{[\theta]}\theta)\cdot\vn\varphi_R dx,
\end{eqnarray*}
so we obtain the identity 
\begin{equation}\label{Identite3}
\int_{\R}div(\mathbb{A}_{[\theta]} \theta)(\theta \varphi_R)dx=-\frac{1}{2}\int_{\mathbb{R}^2} (\mathbb{A}_{[\theta]}\theta^2)\cdot\vn\varphi_Rdx.
\end{equation}
\end{itemize}
Thus, with the identities (\ref{Identite2}) and (\ref{Identite3}), we can thus rewrite the formula (\ref{Identite1}) in the following manner:
$$\int_{B_R} |(-\Delta)^{\frac{\alpha}{4}} \theta|^2\varphi_R dx= \int_{\R} (-\Delta)^{\frac{\alpha}{4}} \theta\left[\left((-\Delta)^{\frac{\alpha}{4}}\theta\right) \varphi_R-(-\Delta)^{\frac{\alpha}{4}}(\theta\varphi_R)\right]dx+\frac{1}{2}\int_{\R} (\mathbb{A}_{[\theta]}\theta^2)\cdot\vn\varphi_Rdx.$$
We recall now that since $0\leq \varphi_R(x)\leq 1$ and $\varphi_R (x)=1$ if $|x|< \frac{R}{2}$,  we have the estimate
$$\int_{B_{\frac{R}{2}}} |(-\Delta)^{\frac{\alpha}{4}} \theta|^2 dx \leq  \int_{B_R}|(-\Delta)^{\frac{\alpha}{4}} \theta|^2\varphi_R dx,$$ 
so we can write
\begin{eqnarray}
\int_{B_{\frac{R}{2}}} |(-\Delta)^{\frac{\alpha}{4}} \theta|^2dx&\leq & \underbrace{\int_{\R} (-\Delta)^{\frac{\alpha}{4}} \theta\left[\left((-\Delta)^{\frac{\alpha}{4}}\theta\right) \varphi_R-(-\Delta)^{\frac{\alpha}{4}}(\theta\varphi_R)\right]dx}_{(I_1)}\notag\\
&&+\underbrace{\frac{1}{2}\int_{\R} (\mathbb{A}_{[\theta]}\theta^2)\cdot\vn\varphi_Rdx}_{(I_2)},\label{Identite4}
\end{eqnarray}
and we will prove now that we have $\underset{R\to +\infty}{\lim}|I_1|=0$ and $\underset{R\to +\infty}{\lim}|I_2|=0$.
\subsection{Estimates for the term $I_1$}
For the term $I_1$ of the expression (\ref{Identite4}) we write, by the Cauchy-Schwartz inequality
\begin{eqnarray*}
|I_1|&=&\left|\int_{\R}(-\Delta)^{\frac{\alpha}{4}} \theta\left[\left((-\Delta)^{\frac{\alpha}{4}}\theta\right) \varphi_R-(-\Delta)^{\frac{\alpha}{4}}(\theta\varphi_R)\right]dx\right|\\
&\leq&\|(-\Delta)^{\frac{\alpha}{4}} \theta\|_{L^2(\R)}\left\|
\left((-\Delta)^{\frac{\alpha}{4}}\theta\right) \varphi_R-(-\Delta)^{\frac{\alpha}{4}}(\theta\varphi_R)\right\|_{L^2(\R)}\\
&\leq &\|\theta\|_{\dot{H}^{\frac{\alpha}{2}}(\R)}\left(\left\|(-\Delta)^{\frac{\alpha}{4}}\big(\theta\varphi_R\big)-\left((-\Delta)^{\frac{\alpha}{4}}\theta\right) \varphi_R-\left((-\Delta)^{\frac{\alpha}{4}} \varphi_R\right) \theta\right\|_{L^2(\R)}\right.\\
&&\qquad\qquad\qquad +\left.\left\|\big((-\Delta)^{\frac{\alpha}{4}} \varphi_R\big) \theta\right\|_{L^2(\R)}\right).
\end{eqnarray*}
We apply now the second point of Lemma \ref{FracLeibniz} to obtain the estimate
\begin{equation}\label{ApplicationKatoPonce}
|I_1|\leq \|\theta\|_{\dot{H}^{\frac{\alpha}{2}}(\R)}\left(\|(-\Delta)^{\frac{\alpha_1}{4}}\varphi_R\|_{L^{p_1}(\R)}\|(-\Delta)^{\frac{\alpha_2}{4}}  \theta\|_{L^{p_2}(\R)}+\left\|\big((-\Delta)^{\frac{\alpha}{4}} \varphi_R\big) \theta\right\|_{L^2(\R)}\right).
\end{equation}
At this point we use the hypothesis $\theta\in L^{\frac{2n-\epsilon}{n-\alpha}}(\R)$ (where $0<\epsilon\ll 1$), thus by the H\"older inequalities with $\frac{1}{2}=\frac{2\alpha-\epsilon}{4n-2\epsilon}+\frac{n-\alpha}{2n-\epsilon}$, we can write:
$$|I_1|\leq \|\theta\|_{\dot{H}^{\frac{\alpha}{2}}(\R)}\left(\|(-\Delta)^{\frac{\alpha_1}{4}}\varphi_R\|_{L^{p_1}(\R)}\|(-\Delta)^{\frac{\alpha_2}{4}}  \theta\|_{L^{p_2}(\R)}+\left\|(-\Delta)^{\frac{\alpha}{4}} \varphi_R\right\|_{L^{\frac{4n-2\epsilon}{2\alpha-\epsilon}}(\R)}\| \theta\|_{L^{\frac{2n-\epsilon}{n-\alpha}}(\R)}\right),$$
now, by the last point of the Lemma \ref{LemmaN2} we have $\|(-\Delta)^{\frac{\alpha_1}{4}}\varphi_R\|_{L^{p_1}(\R)}\leq CR^{-\frac{\alpha_1}{2}+\frac{n}{p_1}}$ as well as $\|(-\Delta)^{\frac{\alpha}{4}} \varphi_R\|_{L^{\frac{4n-2\epsilon}{2\alpha-\epsilon}}(\R)}\leq CR^{-\frac{\alpha}{2}+n(\frac{2\alpha-\epsilon}{4n-2\epsilon})}$ and thus we obtain
$$|I_1|\leq \|\theta\|_{\dot{H}^{\frac{\alpha}{2}}(\R)}\left(CR^{-\frac{\alpha_1}{2}+\frac{n}{p_1}}\|(-\Delta)^{\frac{\alpha_2}{4}}  \theta\|_{L^{p_2}(\R)}+CR^{-\frac{\alpha}{2}+n(\frac{2\alpha-\epsilon}{4n-2\epsilon})}\| \theta\|_{L^{\frac{2n-\epsilon}{n-\alpha}}(\R)}\right).$$
Let us also note that, due to the complex interpolation theory (see \cite[Theorem 6.4.5.]{BL}) we have 
$$\big[\dot{H}^{\frac{\alpha}{2}}, L^{\frac{2n-\epsilon}{n-\alpha}}\big]_{\nu}=\dot{W}^{\frac{\alpha_2}{2}, p_2}\qquad \mbox{and}\qquad \|(-\Delta)^{\frac{\alpha_2}{4}}  \theta\|_{L^{p_2}(\R)}=\|\theta\|_{\dot{W}^{\frac{\alpha_2}{2},p_2}(\R)} \leq C\|\theta\|_{\dot{H}^{\frac{\alpha}{2}}(\R)}^\nu\|\theta\|_{L^\frac{2n-\epsilon}{n-\alpha}(\R)}^{1-\nu},$$
with the relationships 
\begin{equation}\label{ConditionInterpolN}
\alpha_2=\nu\alpha, \qquad \frac{1}{p_2}=\frac{\nu}{2}+(1-\nu)\frac{n-\alpha}{2n-\epsilon}\qquad \mbox{for some} \qquad 0<\nu<1,
\end{equation}
so we can write
$$|I_1|\leq \|\theta\|_{\dot{H}^{\frac{\alpha}{2}}(\R)}\left(CR^{-\frac{\alpha_1}{2}+\frac{n}{p_1}}\|\theta\|_{\dot{H}^{\frac{\alpha}{2}}(\R)}^\nu\|\theta\|_{L^\frac{2n-\epsilon}{n-\alpha}(\R)}^{1-\nu}+CR^{-\frac{\alpha}{2}+n(\frac{2\alpha-\epsilon}{4n-2\epsilon})}\| \theta\|_{L^{\frac{2n-\epsilon}{n-\alpha}}(\R)}\right).$$
But since we have $\alpha=\alpha_1+\alpha_2$ and $\frac{1}{2}=\frac{1}{p_1}+\frac{1}{p_2}$, following the conditions (\ref{ConditionInterpolN}) above, we obtain that $\alpha_1=(1-\nu)\alpha$ and $\frac{1}{p_1}=(1-\nu)\frac{2\alpha-\epsilon}{4n-2\epsilon}$ and we can write 
$$|I_1|\leq \|\theta\|_{\dot{H}^{\frac{\alpha}{2}}(\R)}\left(CR^{(1-\nu)[\frac{(\alpha-n)\epsilon}{4n-2\epsilon}]}\|\theta\|_{\dot{H}^{\frac{\alpha}{2}}(\R)}^\nu\|\theta\|_{L^\frac{2n-\epsilon}{n-\alpha}(\mathbb{R}^2)}^{1-\nu}+CR^{[\frac{(\alpha-n)\epsilon}{4n-2\epsilon}]}\|\theta\|_{L^{\frac{2n-\epsilon}{n-\alpha}}(\R)}\right).$$
We remark now that $\alpha-n<0$ since $1<\alpha<2$ and $n\geq 2$, thus we obtain that the powers $(1-\nu)[\frac{(\alpha-n)\epsilon}{4n-2\epsilon}]$ and $[\frac{(\alpha-n)\epsilon}{4n-2\epsilon}]$ of the parameter $R$ in the previous expression are negative. Moreover, since $\|\theta\|_{\dot{H}^{\frac{\alpha}{2}}(\R)}<+\infty$ and $\|\theta\|_{L^{\frac{2n-\epsilon}{n-\alpha}}(\R)}<+\infty$ by hypothesis, we can write 
\begin{eqnarray*}
\underset{R\to+\infty}{\lim }|I_1|&\leq &C\|\theta\|_{\dot{H}^{\frac{\alpha}{2}}(\R)}\left[\|\theta\|_{\dot{H}^{\frac{\alpha}{2}}(\R)}^\nu\|\theta\|_{L^\frac{2n-\epsilon}{n-\alpha}(\R)}^{1-\nu}\left( \underset{R\to+\infty}{\lim }R^{(1-\nu)[\frac{(\alpha-n)\epsilon}{4n-2\epsilon}]}\right)\right.\\[2mm]
&&\left.+\|\theta\|_{L^{\frac{2n-\epsilon}{n-\alpha}}(\R)}\left( \underset{R\to+\infty}{\lim }R^{[\frac{(\alpha-n)\epsilon}{4n-2\epsilon}]}\right)\right]=0,
\end{eqnarray*}
from which we easily deduce that 
\begin{equation}\label{LimiteI1}
\underset{R\to+\infty}{\lim }|I_1|=0.
\end{equation}
\begin{Remarque}
Note that the study of the term $I_1$ of (\ref{Identite4}) is independent from the dimension $n\geq 2$ and it is valid for all $0<\alpha<2$. Note also that the hypothesis $\theta\in L^{\frac{2n-\epsilon}{n-\alpha}}(\R)$ with a small $0<\epsilon \ll1$ is only useful to apply the Kato-Ponce inequalities in the formula (\ref{ApplicationKatoPonce}) above and to obtain some interpolation estimates that provide the suitable information in order to pass to the limit $R\to+\infty$.
\end{Remarque}
\subsection{Estimates for the term $I_2$ of (\ref{Identite4})}
We need to estimate the quantity 
$$|I_2|\leq \int_{\R} \left|(\mathbb{A}_{[\theta]}\theta^2)\cdot\vn\varphi_R\right|dx,$$
and for this we will distinguish some cases following the values of the parameter $0<\alpha<2$ and the dimension $n\geq 2$. Indeed: 
\begin{itemize}
\item[$\bullet$] If $\frac{n}{3}<\alpha<\frac{n+2}{3}\leq 2$, ($n=2,3,4$), by the H\"older inequality with 
\begin{equation}\label{IneqHolder1}
1=\frac{n-\alpha}{2n}+\frac{n-\alpha}{2n}+\frac{n-\alpha}{2n}+\frac{3\alpha-n}{2n},
\end{equation}
we can write 
\begin{eqnarray*}
|I_2|&\leq &\|\mathbb{A}_{[\theta]}\|_{L^{\frac{2n}{n-\alpha}}(\R)}\|\theta\|_{L^{\frac{2n}{n-\alpha}}(\R)}^2 \|\vn\varphi_R\|_{L^{\frac{2n}{3\alpha-n}}(\R)}\\
&\leq & C_{\mathbb{A}}\|\theta\|_{L^{\frac{2n}{n-\alpha}}(\R)} \|\theta\|_{L^{\frac{2n}{n-\alpha}}(\R)}^2 \|\vn\varphi_R\|_{L^{\frac{2n}{3\alpha-n}}(\R)},
\end{eqnarray*}
where in the last line we used the boundedness property of the operator $\mathbb{A}$ in Lebesgue spaces (\emph{i.e.}. $\|\mathbb{A}_{[\theta]}\|_{L^{\frac{2n}{n-\alpha}}(\R)}\leq C_{\mathbb{A}}\|\theta\|_{L^{\frac{2n}{n-\alpha}}(\R)}$). Now, by the point 4) of the Lemma \ref{LemmaN2}, we have $\|\vn\varphi_R\|_{L^{\frac{2n}{3\alpha-n}}(\R)}\leq CR^{-1+n(\frac{3\alpha-n}{2n})}$ and we can write 
$$|I_2|\leq C_{\mathbb{A}}\|\theta\|_{L^{\frac{2n}{n-\alpha}}(\R)}^3 R^{-1+n(\frac{3\alpha-n}{2n})}.$$
We note now that since $\alpha<\frac{n+2}{3}$ we have $-1+n(\frac{3\alpha-n}{2n})<0$ moreover recalling that by the Sobolev inequalities we have $\|\theta\|_{L^{\frac{2n}{n-\alpha}}(\R)}\leq C\|\theta\|_{\dot{H}^{\frac{\alpha}{2}}(\R)}<+\infty$, then we obtain 
$$\underset{R\to+\infty}{\lim }|I_2|=0.$$
\item[$\bullet$] If $\alpha=\frac{n}{3}<2$, ($n=2,3,4$), then we proceed as follows: recalling that we have the embedding $\dot{H}^{\frac{n}{6}}(\R)\subset L^3(\R)$ we write by the H\"older inequalities with 
\begin{equation}\label{IneqHolder2}
1=\frac{1}{3}+\frac{1}{3}+\frac{1}{3},
\end{equation}
\begin{eqnarray*}
|I_2|&\leq &\int_{\R} |\mathbb{A}_{[\theta]}||\theta|^2|\vn\varphi_R|dx\leq \|\mathbb{A}_{[\theta]}\|_{L^3(\R)}\|\theta\|_{L^3(\R)}^2\|\vn\varphi_R\|_{L^\infty(\R)}\\
&\leq &C_{\mathbb{A}}\|\theta\|_{L^3(\R)}^3 CR^{-1},
\end{eqnarray*}
where we used the boundedness property of the operator $\mathbb{A}$ in Lebesgue spaces as well as the properties of the test function $\varphi_R$ given in the point 4) of the Lemma \ref{LemmaN2}. From this estimate, since we have $\|\theta\|_{L^3(\R)}\leq \|\theta\|_{\dot{H}^{\frac{n}{6}}(\R)}<+\infty$, we easily deduce that 
$$\underset{R\to+\infty}{\lim }|I_2|=0.$$
\item[$\bullet$] If $\alpha=\frac{n+2}{3}<2$, ($n=2,3$), then we have $\dot{H}^{\frac{n+2}{6}}(\R)\subset L^{\frac{3n}{n-1}}(\R)$ and recalling that $supp(\nabla \varphi_R)\subset \mathcal{C}_R = \{ x \in \mathbb{R}^n : \tfrac{R}{2} \leq |x| \leq R \}$ we write, by the H\"older inequalities with 
\begin{equation}\label{IneqHolder3}
1=\frac{n-1}{3n}+\frac{n-1}{3n}+\frac{n-1}{3n}+\frac{1}{n},
\end{equation} 
\begin{eqnarray*}
|I_2|&\leq &\int_{\R} |\mathbb{A}_{[\theta]}||\theta|^2\mathds{1}_{\mathcal{C}_R}|\vn\varphi_R|dx\leq \|\mathbb{A}_{[\theta]}\|_{L^{\frac{3n}{n-1}}(\R)}\|\theta\|_{L^\frac{3n}{n-1}(\mathcal{C}_R)}^2\|\vn\varphi_R\|_{L^n(\R)}\\
&\leq &C_{\mathbb{A}}\|\theta\|_{L^{\frac{3n}{n-1}}(\R)}\|\theta\|_{L^{\frac{3n}{n-1}}(\mathcal{C}_R)}^2,
\end{eqnarray*}
since in this case we have $\|\vn\varphi_R\|_{L^n(\R)}\leq C$. Now, as we have by the Sobolev embeddings the uniform control $\|\theta\|_{L^{\frac{3n}{n-1}}(\R)}\leq C\|\theta\|_{\dot{H}^{\frac{n+2}{6}}(\R)}<+\infty$ we obtain that  $\underset{R\to+\infty}{\lim }\|\theta\|_{L^{\frac{3n}{n-1}}(\mathcal{C}_R)}=0$, from which we deduce that 
$$\underset{R\to+\infty}{\lim }|I_2|=0.$$

\begin{Remarque}
Note that in all the previous cases the ``natural'' Lebesgue information $L^{\frac{2n}{n-\alpha}}(\R)$ (obtained by the Sobolev embedding $\dot{H}^{\frac{\alpha}{2}}(\R)\subset L^{\frac{2n}{n-\alpha}}(\R)$) provides enough information to obtain $\underset{R\to+\infty}{\lim }|I_2|=0$. In order to deduce this limit outside these cases we will ask some additional information over the solution $\theta \in \dot{H}^{\frac{\alpha}{2}}(\R)$.
\end{Remarque}

\item[$\bullet$] If $0<\alpha<\frac{n}{3}<2$, ($n=2,3,4$), in this case the Lebesgue information given by the Sobolev embedding $\dot{H}^{\frac{\alpha}{2}}(\R)\subset L^{\frac{2n}{n-\alpha}}(\R)$ is not well suited to use the H\"older inequalities as in (\ref{IneqHolder1})-(\ref{IneqHolder3}) above, and we will thus ask the additional condition $\theta \in L^{\frac{2n+\nu}{n-\alpha}}(\R)$ for some real parameter $\nu$ such that  $(n-3\alpha)<\nu <1+(n-3\alpha)$. Then we have, by the H\"older inequality with $1=\frac{n-\alpha}{2n+\nu}+\frac{n-\alpha}{2n+\nu}+\frac{n-\alpha}{2n+\nu}+\frac{3\alpha+\nu-n}{2n+\nu}$: 
\begin{eqnarray*}
|I_2|&\leq &\int_{\R} |\mathbb{A}_{[\theta]}||\theta|^2|\vn\varphi_R|dx\leq \|\mathbb{A}_{[\theta]}\|_{L^{\frac{2n+\nu}{n-\alpha}}(\R)}\|\theta\|_{L^{\frac{2n+\nu}{n-\alpha}}(\R)}^2\|\vn\varphi_R\|_{L^{\frac{2n+\nu}{3\alpha+\nu-n}}(\R)}\\
&\leq &C_{\mathbb{A}}\|\theta\|_{L^{\frac{2n+\nu}{n-\alpha}}(\R)}^3 CR^{-1+n(\frac{3\alpha+\nu-n}{2n+\nu})},
\end{eqnarray*}
but since $n(\frac{3\alpha+\nu-n}{2n+\nu})<1$ as we have $\nu <1+(n-3\alpha)$, we easily obtain
$$\underset{R\to+\infty}{\lim }|I_2|=0.$$

\item[$\bullet$] If $\frac{n+2}{3}\leq \alpha<2$, ($n=2,3$), the Sobolev embeddings do not provide enough information and we will ask $\theta \in L^{\frac{3n}{n-1}}(\R)$, and thus, following the same ideas as in (\ref{IneqHolder3}) above, we obtain $\underset{R\to+\infty}{\lim }|I_2|=0$.
\begin{Remarque}
Note that in this last case, in dimension $n=3$, the additional hypotheses asked here is $L^{\frac{9}{2}}(\R)$ which corresponds to the best Lebesgue space known to date when we consider the Liouville-type theorems for the Navier-Stokes equations, see \cite[Theorem X.9.5]{Galdi}.
\end{Remarque}
\end{itemize}
 \hfill $\blacksquare$
\appendix
\mysection{Weak solutions in $\dot{H}^{\frac{\alpha}{2}}(\R)$}\label{Secc_ProofTheo1}

In order to solve equation (\ref{SQGStationnaire}) we will first consider a function $\phi\in \mathcal{C}^{\infty}_{0}(\R)$ such that $0\leq \phi (x)\leq 1$ with $\phi(x)=1$ if $|x|\leq 1$ and $\phi (x)=0$ if $|x|>2$, then for $R>1$ we set $\phi_{R}(x)=\phi(\frac{x}{R})$. We also consider a positive test function $\varphi\in \mathcal{C}^{\infty}_{0}(\R)$ and we set $\varphi_\epsilon(x)=\frac{1}{\epsilon^{n}}\varphi(\frac{x}{\epsilon})$ for $\epsilon >0$. 
 With these auxiliary functions and for some $0<\rho<1$ we study the following equation
\begin{equation}\label{EquationRegulariseeStationnaire}
\rho(-\Delta) \theta +(-\Delta)^{\frac \alpha 2}\theta-(\mathbb{A}_{[(\theta \phi_R)\ast \varphi_\epsilon]}\ast \varphi_\epsilon)\phi_R\cdot \vn (\theta \phi_R)-f=0.
\end{equation}
Remark that, at least formally, if we make $\rho, \epsilon\to 0$ and $R\to +\infty$, we recover the equation (\ref{SQGStationnaire}).\\

We note now that from the equation (\ref{EquationRegulariseeStationnaire}) we can write
\begin{equation}\label{ProblemePointFixe}
\theta=\mathbb{T}_{R, \rho, \epsilon}(\theta),
\end{equation}
where
\begin{equation}\label{DefTREpsi}
\mathbb{T}_{R, \rho, \epsilon}(\theta)=\mathfrak{L}_{\rho}\left((\mathbb{A}_{[(\theta \phi_R)\ast \varphi_\epsilon]}\ast \varphi_\epsilon)\phi_R\cdot \vn (\theta \phi_R) +f\right),
\end{equation}
and where the operator $\mathfrak{L}_{\rho}$ is defined, for a fixed $0<\rho<1$ by the expression
\begin{equation}\label{DefopL}
\mathfrak{L}_{\rho}(f)=\frac{1}{[\rho(-\Delta)+(-\Delta)^{\frac{\alpha}{2}}]}(f).
\end{equation}
The operator $\mathfrak{L}_{\rho}$ satisfies the following property.
\begin{Lemme}\label{LemmeControlOpFourier}
Fix $0<\rho<1$ and consider a regularity index $0<\alpha<2$. Then, for all parameter $\sigma$ such that $\alpha\leq \sigma\leq 2$ and for any regular function $f:\R\longrightarrow \mathbb{R}$, we have the estimate
\begin{equation}\label{ControlOpFourier}
\|\mathfrak{L}_{\rho}(f)\|_{\dot{H}^{\sigma}(\R)}\leq \frac{C}{\rho}\|f\|_{L^2(\R)}.
\end{equation}
\end{Lemme}
Indeed, since we are working in the $L^2$ space, the proof of this estimate follows from the Plancherel theorem and the definition of the homogeneous Sobolev space $\dot{H}^\sigma$. We can thus write
$$\|\mathfrak{L}_{\rho}(f)\|_{\dot{H}^{\sigma}(\R)}=\left\|\frac{(-\Delta)^{\frac{\sigma}{2}}}{[\rho(-\Delta)+(-\Delta)^{\frac{\alpha}{2}}]}(f)\right\|_{L^2(\R)}=\left\|\frac{|\xi|^{\sigma}}{[\rho|\xi|^2+|\xi|^{\alpha}]}\widehat{f}(\cdot)\right\|_{L^2(\R)}.$$
We remark now that, since $\alpha\leq \sigma\leq 2$, we have the uniform estimate $\frac{|\xi|^{\sigma}}{[\rho|\xi|^2+|\xi|^{\alpha}]}\leq \frac{C}{\rho}$, and we obtain 
$$\|\mathfrak{L}_{\rho}(f)\|_{\dot{H}^{\sigma}(\R)}\leq \frac{C}{\rho}\left\|\widehat{f}(\cdot)\right\|_{L^2(\R)}=\frac{C}{\rho}\|f\|_{L^2(\R)},$$
which is the wished estimate.\\

Now, in order to obtain a solution for the problem $\theta=\mathbb{T}_{R, \rho}(\theta)$ we will apply the Schaefer fixed point theorem (see \cite[Theorem 16.1]{PGLR2}):
\begin{Theoreme}[Schaefer]\label{Teo_Schaefer}
For $1<\alpha<2$, consider the following functional space:
\begin{equation}\label{DefFuncSpace}
E_{\rho}=\big\{\theta:\R\longrightarrow \mathbb{R}: \theta\in \dot{H}^{1}(\R)\cap\dot{H}^{\frac{\alpha}{2}}(\R) \big\},
\end{equation}
endowed with the norm 
\begin{equation}\label{DefNormeE}
\|\theta\|_{E_{\rho}}=\sqrt{\rho}\|\theta\|_{\dot{H}^{1}(\R)}+\|\theta\|_{\dot{H}^{\frac{\alpha}{2}}(\R)}.
\end{equation}
If we have the following points:
\begin{itemize}
\item[$1)$] the operator $\mathbb{T}_{R, \rho, \epsilon}$ defined in (\ref{DefTREpsi}) is continuous and compact in the space $E_{\rho}$,
\item[$2)$] if $\theta=\lambda \mathbb{T}_{R, \rho, \epsilon}(\theta)$ for any $\lambda\in [0,1]$, then we have $\|\theta\|_{E_{\rho}}\leq M$,
\end{itemize}
then the equation (\ref{ProblemePointFixe}) admits at least one solution $\theta\in E_{\rho}$.
\end{Theoreme}
As we can see, to obtain a solution of the modified problem (\ref{EquationRegulariseeStationnaire}), it is enough to verify the two points of the previous theorem. We decompose our study in some propositions and corollaries that will be helpful in the sequel. 

\begin{Proposition}[Continuity]\label{PropoContinuidad}
The application $\mathbb{T}_{R, \rho, \epsilon}$ is continuous in the space $E_{\rho}$.
\end{Proposition}
{\bf Proof.} We start writing
\begin{eqnarray}
\|\mathbb{T}_{R, \rho, \epsilon}(\theta)\|_{E_{\rho}}&=&\sqrt{\rho}\|\mathbb{T}_{R, \rho, \epsilon}(\theta)\|_{\dot{H}^1(\R)}+\|\mathbb{T}_{R, \rho, \epsilon}(\theta)\|_{\dot{H}^{\frac{\alpha}{2}}(\R)}\notag\\
&=&\sqrt{\rho}\underbrace{\left\|\mathfrak{L}_{\rho}\left((\mathbb{A}_{[(\theta \phi_R)\ast \varphi_\epsilon]}\ast \varphi_{\epsilon})\phi_R\cdot \vn (\theta \phi_R) +f\right)\right\|_{\dot{H}^1(\R)}}_{(1)}\notag\\
&&+\underbrace{\left\|\mathfrak{L}_{\rho}\left((\mathbb{A}_{[(\theta \phi_R)\ast \varphi_\epsilon]}\ast \varphi_\epsilon)\phi_R\cdot \vn (\theta \phi_R) +f\right)\right\|_{\dot{H}^{\frac{\alpha}{2}}(\R)}}_{(2)},\label{Continuite0}
\end{eqnarray}
and we will study each term above separately. 
\begin{itemize}
\item[$\bullet$] We will bound the term (1) above by the quantity $\|\theta\|_{E_\rho}$ and the information available over the external force $f$. We thus write
$$(1)\leq \left\|\mathfrak{L}_{\rho}\left((\mathbb{A}_{[\theta \phi_R]}\ast \varphi_{\epsilon})\phi_R\cdot \vn (\theta \phi_R) \right)\right\|_{\dot{H}^1(\R)}+\left\|\mathfrak{L}_{\rho}\left(f\right)\right\|_{\dot{H}^1(\R)},$$
and by the definition of the operator $\mathfrak{L}_{\rho}$ given in (\ref{DefopL}), we write
\begin{eqnarray*}
(1)\leq \left\| \frac{1}{[\rho(-\Delta)+(-\Delta)^{\frac{\alpha}{2}}]}\left((\mathbb{A}_{[(\theta \phi_R)\ast \varphi_\epsilon]}\ast \varphi_{\epsilon})\phi_R\cdot \vn (\theta \phi_R) \right)\right\|_{\dot{H}^1(\R)}+\left\|  \frac{1}{[\rho(-\Delta)+(-\Delta)^{\frac{\alpha}{2}}]}\left(f\right)\right\|_{\dot{H}^1(\R)},
\end{eqnarray*}
now for some parameter $\sigma>0$ such that  $0<\frac{2-\alpha}{2}<\sigma<1$ and $\alpha<1+\sigma<2$ and by the definition of the Sobolev space $\dot{H}^1(\R)$, we have
\begin{eqnarray*}
(1)\leq \left\| \frac{(-\Delta)^{\frac{1+\sigma}{2}}}{[\rho(-\Delta)+(-\Delta)^{\frac{\alpha}{2}}]}\frac{1}{(-\Delta)^{\frac{\sigma}{2}}}\left((\mathbb{A}_{[(\theta \phi_R)\ast \varphi_\epsilon]}\ast \varphi_{\epsilon})\phi_R\cdot \vn (\theta \phi_R) \right)\right\|_{L^2(\R)}\\
+\left\|  \frac{(-\Delta)^{\frac{1+\frac{\alpha}{2}}{2}}}{[\rho(-\Delta)+(-\Delta)^{\frac{\alpha}{2}}]}\frac{1}{(-\Delta)^{\frac{\alpha}{4}}}\left(f\right)\right\|_{L^2(\R)}.
\end{eqnarray*}
We can thus apply the Lemma \ref{LemmeControlOpFourier} (since we have $\alpha<1+\sigma<2$ and $\alpha<1+\frac{\alpha}{2}<2$) to obtain the estimate
\begin{eqnarray*}
(1)&\leq &\frac{C}{\rho}\left\|\frac{1}{(-\Delta)^{\frac{\sigma}{2}}}\left((\mathbb{A}_{[(\theta \phi_R)\ast \varphi_\epsilon]}\ast \varphi_{\epsilon})\phi_R\cdot \vn (\theta \phi_R) \right)\right\|_{L^2(\R)}+\frac{C}{\rho}\left\|\frac{1}{(-\Delta)^{\frac{\alpha}{4}}}\left(f\right)\right\|_{L^2(\R)}\\
&\leq &\frac{C}{\rho}\left\|(\mathbb{A}_{[(\theta \phi_R)\ast \varphi_\epsilon]}\ast \varphi_{\epsilon})\phi_R\cdot \vn (\theta \phi_R) \right\|_{\dot{H}^{-\sigma}(\R)}+\frac{C}{\rho}\left\|f\right\|_{\dot{H}^{-\frac{\alpha}{2}}(\R)},
\end{eqnarray*}
where in the last control we used the definition of the Sobolev space $\dot{H}^{-\sigma}(\R)$. At this point we apply the Hardy-Littlewood-Sobolev embedding $L^{\frac{2n}{n+2\sigma}}(\R)\subset \dot{H}^{-\sigma}(\R)$ (recall that $0<\sigma<1$) and we have
$$(1)\leq \frac{C}{\rho}\left\|(\mathbb{A}_{[(\theta \phi_R)\ast \varphi_\epsilon]}\ast \varphi_{\epsilon})\phi_R\cdot \vn (\theta \phi_R) \right\|_{L^{\frac{2n}{n+2\sigma}}(\R)}+\frac{C}{\rho}\left\|f\right\|_{\dot{H}^{-\frac{\alpha}{2}}(\R)},$$
and by the H\"older inequalities with $\frac{n+2\sigma}{2n}=\frac{\sigma}{n}+\frac{1}{2}$ we can thus write
\begin{eqnarray*}
(1)&\leq &\frac{C}{\rho}\left\|(\mathbb{A}_{[(\theta \phi_R)\ast \varphi_\epsilon]}\ast \varphi_{\epsilon})\phi_R\right\|_{L^{\frac{n}{\sigma}}(\R)} \left\|\vn (\theta \phi_R) \right\|_{L^{2}(\R)}+\frac{C}{\rho}\left\|f\right\|_{\dot{H}^{-\frac{\alpha}{2}}(\R)}\\
&\leq &\frac{C}{\rho}\|\phi_R\|_{L^\infty(\R)}\left\|\mathbb{A}_{[(\theta \phi_R)\ast \varphi_\epsilon]}\ast \varphi_{\epsilon}\right\|_{L^{\frac{n}{\sigma}}(\R)} \left(\|(\vn \theta) \phi_R \|_{L^{2}(\R)}+\|(\vn \phi_R) \theta \|_{L^{2}}\right)+\frac{C}{\rho}\left\|f\right\|_{\dot{H}^{-\frac{\alpha}{2}}(\R)}.
\end{eqnarray*}
We apply now the Young inequalities with $1+\frac{\sigma}{n}=\frac{1}{n}+\frac{n+\sigma-1}{n}$ and the  H\"older inequalities with $\frac{1}{2}=\frac{n-\alpha}{2n}+\frac{\alpha}{2n}$ to obtain
\begin{eqnarray*}
(1)&\leq& \frac{C}{\rho}\left\|\mathbb{A}_{[(\theta \phi_R)\ast \varphi_\epsilon]}\right\|_{L^{n}(\R)}\|\varphi_\epsilon\|_{L^{\frac{n}{n+\sigma-1}}(\R)} \left(\|\phi_R\|_{L^\infty}\|\vn \theta  \|_{L^{2}(\R)}+\|\vn \phi_R\|_{L^\frac{2n}{\alpha}(\R)}\|\theta \|_{L^{\frac{2n}{n-\alpha}}(\R)}\right)\\
&&+\frac{C}{\rho}\left\|f\right\|_{\dot{H}^{-\frac{\alpha}{2}}(\R)}\\
&\leq & \frac{C_{R,\epsilon}}{\rho}\left\|(\theta \phi_R)\ast \varphi_\epsilon\right\|_{L^{n}(\R)}\left(\|\vn\theta\|_{L^{2}(\R)}+\|\theta \|_{L^{\frac{2n}{n-\alpha}}(\R)}\right)+\frac{C}{\rho}\left\|f\right\|_{\dot{H}^{-\frac{\alpha}{2}}(\R)},
\end{eqnarray*}
where we used the fact that the operator $\mathbb{A}$ is bounded in the Lebesgue spaces.\\

Now, by the Young inequalities we have 
\begin{eqnarray*}
(1)&\leq & \frac{C_{R,\epsilon}}{\rho}\|\theta \phi_R\|_{L^1(\R)}\|\varphi_\epsilon\|_{L^{n}(\R)}\left(\|\vn\theta\|_{L^{2}(\R)}+\|\theta \|_{L^{\frac{2n}{n-\alpha}}(\R)}\right)+\frac{C}{\rho}\left\|f\right\|_{\dot{H}^{-\frac{\alpha}{2}}(\R)}\\
&\leq &  \frac{C_{R,\epsilon}}{\rho}\|\theta\|_{L^{\frac{2n}{n-\alpha}}(\R)}\|\phi_R\|_{L^{\frac{2n}{n+\alpha}}(\R)}\left(\|\vn\theta\|_{L^{2}(\R)}+\|\theta \|_{L^{\frac{2n}{n-\alpha}}(\R)}\right)+\frac{C}{\rho}\left\|f\right\|_{\dot{H}^{-\frac{\alpha}{2}}(\R)},
\end{eqnarray*}
where we used the H\"older inequalities with $1=\frac{n-\alpha}{2n}+\frac{n+\alpha}{2n}$ (recall that $0<\alpha<2$ and that $n\geq 2$). Now, by the Sobolev embedding $\dot{H}^{\frac{\alpha}{2}}(\R)\subset L^{\frac{2n}{n-\alpha}}(\R)$, we can write
\begin{eqnarray*}
(1)&\leq &\frac{C_{R,\epsilon}}{\rho}\|\theta \|_{\dot{H}^{\frac{\alpha}{2}}(\R)}\left(\|\theta\|_{\dot{H}^{1}(\R)}+\|\theta \|_{\dot{H}^{\frac{\alpha}{2}}(\R)}\right)+\frac{C}{\rho}\left\|f\right\|_{\dot{H}^{-\frac{\alpha}{2}}(\R)}.
\end{eqnarray*}
We remark that we have $\|\cdot\|_{\dot{H}^{\frac{\alpha}{2}}(\R)}\leq \|\cdot\|_{E_\rho}$ and $\|\cdot\|_{\dot{H}^{1}(\R)}\leq \frac{1}{\sqrt{\rho}}\|\cdot\|_{E_\rho}$ (see formula (\ref{DefNormeE})), so we can write 
\begin{eqnarray}
(1)\leq \frac{C_{R,\epsilon}}{\rho}\|\theta \|_{E_\rho}\left(\frac{1}{\sqrt{\rho}}\|\theta\|_{E_\rho}+\|\theta \|_{E_\rho}\right)+\frac{C}{\rho}\left\|f\right\|_{\dot{H}^{-\frac{\alpha}{2}}(\R)}\notag\\
\leq \frac{C_{R,\epsilon}}{\rho}\left(1+\frac{1}{\sqrt{\rho}}\right)\|\theta\|_{E_\rho}\|\theta \|_{E_\rho}+\frac{C}{\rho}\left\|f\right\|_{\dot{H}^{-\frac{\alpha}{2}}(\R)},\label{Estimate1Continuite}
\end{eqnarray}
and this estimate ends the study of the first term of (\ref{Continuite0}).

\item[$\bullet$] For the term (2) in (\ref{Continuite0}) we write
$$(2)\leq \left\|\mathfrak{L}_{\rho}\left((\mathbb{A}_{[(\theta \phi_R)\ast \varphi_\epsilon]}\ast \varphi_{\epsilon})\phi_R\cdot \vn (\theta \phi_R)\right)\right\|_{\dot{H}^{\frac{\alpha}{2}}(\R)}+\left\|\mathfrak{L}_{\rho}\left(f\right)\right\|_{\dot{H}^{\frac{\alpha}{2}}(\R)}.$$

Using the definition of the Sobolev spaces $\dot{H}^{\frac{\alpha}{2}}(\R)$ and the definition of the operator $\mathfrak{L}_{\rho}$ given in (\ref{DefopL}), for some $\frac{\alpha}{2}<s<1$ (we thus have $\frac{\alpha}{2}<\frac{s+\frac{\alpha}{2}}{2}<2$), we can write
\begin{eqnarray*}
(2)&\leq &\left\| \frac{(-\Delta)^{\frac{s+\frac{\alpha}{2}}{2}}}{[\rho(-\Delta)+(-\Delta)^{\frac{\alpha}{2}}]}\frac{1}{(-\Delta)^{\frac{s}{2}}}\left((\mathbb{A}_{[(\theta \phi_R)\ast \varphi_\epsilon]}\ast \varphi_\epsilon)\phi_R\cdot \vn (\theta \phi_R) \right)\right\|_{L^2(\R)}\\
&&+\left\|\frac{(-\Delta)^{\frac{\alpha}{2}}}{[\rho(-\Delta)+(-\Delta)^{\frac{\alpha}{2}}]}\frac{1}{(-\Delta)^{\frac{\alpha}{4}}}\left(f\right)\right\|_{L^2(\R)},
\end{eqnarray*}
and since $\frac{\alpha}{2}<\frac{s+\frac{\alpha}{2}}{2}<2$ we can apply the Lemma \ref{LemmeControlOpFourier} to obain
\begin{eqnarray*}
(2)&\leq &\frac{C}{\rho}\left\|\frac{1}{(-\Delta)^{\frac{s}{2}}}\left((\mathbb{A}_{[(\theta \phi_R)\ast \varphi_\epsilon]}\ast \varphi_\epsilon)\phi_R\cdot \vn (\theta \phi_R) \right)\right\|_{L^2(\R)}+\frac{C}{\rho}\left\|\frac{1}{(-\Delta)^{\frac{\alpha}{4}}}\left(f\right)\right\|_{L^2(\R)}\\
&\leq &\frac{C}{\rho}\left\|(\mathbb{A}_{[(\theta \phi_R)\ast \varphi_\epsilon]}\ast \varphi_\epsilon)\phi_R\cdot \vn (\theta \phi_R) \right\|_{\dot{H}^{-s}(\R)}+\frac{C}{\rho}\left\|f\right\|_{\dot{H}^{-\frac{\alpha}{2}}(\R)}.
\end{eqnarray*}

By the Hardy-Littlewood-Sobolev embedding $L^{\frac{2n}{n+2s}}(\R)\subset \dot{H}^{-s}(\R)$ (recall that $\frac{\alpha}{2}<s<1\leq \frac{n}{2}$ for $n\geq 2$), we obtain
\begin{eqnarray*}
(2)&\leq &\frac{C}{\rho}\left\|(\mathbb{A}_{[(\theta \phi_R)\ast \varphi_\epsilon]}\ast \varphi_\epsilon)\phi_R\cdot \vn (\theta \phi_R) \right\|_{L^{\frac{2n}{n+2s}}(\R)}+\frac{C}{\rho}\left\|f\right\|_{\dot{H}^{-\frac{\alpha}{2}}(\R)},
\end{eqnarray*}
now, by the H\"older inequality with $\frac{n+2s}{2n}=\frac{s}{n}+\frac{1}{2}$, we write
\begin{eqnarray*}
(2)&\leq &\frac{C}{\rho}\left\|(\mathbb{A}_{[(\theta \phi_R)\ast \varphi_\epsilon]}\ast \varphi_\epsilon)\phi_R\right\|_{L^\frac{n}{s}(\R)} \|\vn (\theta \phi_R)\|_{L^{2}(\R)}+\frac{C}{\rho}\left\|f\right\|_{\dot{H}^{-\frac{\alpha}{2}}(\R)}\\
&\leq &\frac{C}{\rho}\|\phi_R\|_{L^\infty(\R)}\|\mathbb{A}_{[(\theta \phi_R)\ast \varphi_\epsilon]}\|_{L^n(\R)}\|\varphi_\epsilon\|_{L^{\frac{n}{n+s-1}}(\R)} \|\vn (\theta \phi_R)\|_{L^{2}(\R)}+\frac{C}{\rho}\left\|f\right\|_{\dot{H}^{-\frac{\alpha}{2}}(\R)}\\
&\leq &\frac{C_{\epsilon}}{\rho}\|(\theta \phi_R)\ast \varphi_\epsilon\|_{L^n(\R)}\|\vn (\theta \phi_R)\|_{L^{2}(\R)}+\frac{C}{\rho}\left\|f\right\|_{\dot{H}^{-\frac{\alpha}{2}}(\R)}
\end{eqnarray*}
where we used the Young inequalities with $1+\frac{s}{n}=\frac{1}{n}+\frac{n+s-1}{n}$ (recall that $0<s<1$) and the fact that the operator $\mathbb{A}$ is bounded in Lebesgue spaces. Using again the Young inequalities we have 
\begin{eqnarray*}
(2)&\leq &\frac{C_{\epsilon}}{\rho}\|\theta \phi_R\|_{L^1(\R)}\|\varphi_\epsilon\|_{L^n(\R)}\left(\|\phi_R\|_{L^\infty}\|\vn \theta  \|_{L^{2}(\R)}+\|\vn \phi_R\|_{L^\frac{2n}{\alpha}(\R)}\|\theta \|_{L^{\frac{2n}{n-\alpha}}(\R)}\right)\\
&&+\frac{C}{\rho}\left\|f\right\|_{\dot{H}^{-\frac{\alpha}{2}}(\R)}\\
&\leq &\frac{C_{R,\epsilon}}{\rho}\|\theta\|_{L^{\frac{2n}{n-\alpha}}(\R)}\|\phi_R\|_{L^{\frac{2n}{n+\alpha}}(\R)}\left(\|\vn \theta  \|_{L^{2}(\R)}+\|\theta \|_{L^{\frac{2n}{n-\alpha}}(\R)}\right)+\frac{C}{\rho}\left\|f\right\|_{\dot{H}^{-\frac{\alpha}{2}}(\R)},
\end{eqnarray*}
where we applied the H\"older inequalities with $1=\frac{n-\alpha}{2n}+\frac{n+\alpha}{2n}$. Now, by the Sobolev embedding $\dot{H}^{\frac{\alpha}{2}}(\R)\subset L^{\frac{2n}{n-\alpha}}(\R)$ we obtain
$$(2)\leq \frac{C_{R,\epsilon}}{\rho}\|\theta\|_{\dot{H}^{\frac{\alpha}{2}}(\R)}\left(\|\vn \theta  \|_{L^{2}(\R)}+\|\theta \|_{\dot{H}^{\frac{\alpha}{2}}(\R)}\right)+\frac{C}{\rho}\left\|f\right\|_{\dot{H}^{-\frac{\alpha}{2}}(\R)}.$$
Recalling that $\|\cdot\|_{\dot{H}^{\frac{\alpha}{2}}(\R)}\leq \|\cdot\|_{E_\rho}$ and $\|\cdot\|_{\dot{H}^{1}(\R)}\leq \frac{1}{\sqrt{\rho}}\|\cdot\|_{E_\rho}$, we obtain
\begin{equation}\label{Estimate2Continuite}
(2)\leq\frac{C_R}{\rho}\left(1+\frac{1}{\sqrt{\rho}}\right)\|\theta\|_{E_\rho}\|\theta \|_{E_\rho}+\frac{C}{\rho}\left\|f\right\|_{\dot{H}^{-\frac{\alpha}{2}}(\R)}.
\end{equation}
\end{itemize}
Now, gathering the estimates (\ref{Estimate1Continuite}) and (\ref{Estimate2Continuite}) and coming back to (\ref{Continuite0}), we have proven the control
$$\|\mathbb{T}_{R, \rho, \epsilon}(\theta)\|_{E_{\rho}}\leq \frac{C_R}{\sqrt{\rho}}\left(1+\frac{1}{\sqrt{\rho}}\right)\|\theta\|_{E_\rho}\|\theta \|_{E_\rho}+\frac{C}{\sqrt{\rho}}\left\|f\right\|_{\dot{H}^{-\frac{\alpha}{2}}(\R)},$$
which gives the continuity of the operator $\mathbb{T}_{R, \rho, \epsilon}$ in the space $E_\rho$ (recall that by hypothesis we have $f\in \dot{H}^{-\frac{\alpha}{2}}(\R)$) and Proposition \ref{PropoContinuidad} is now proven. \hfill $\blacksquare$

\begin{Proposition}[Compactness]\label{PropoCompacidad}
The application $\mathbb{T}_{R, \rho, \epsilon}$ is compact in the space $E_{\rho}$.
\end{Proposition}
{\bf Proof.} We consider now a bounded sequence $(\theta_n )_{n\in \mathbb{N}}$ in $E_\rho$ (\emph{i.e.} we have $\|\theta_n\|_{E_\rho}<+\infty$ uniformly) and we shall prove that there exists a subsequence $\mathbb{T}_{R, \rho, \epsilon}(\theta_{n_k})_{k\in \mathbb{N}}$ which converges strongly in $E_\rho$.\\

We first remark that for some fixed $R>0$, the sequence $(\theta_n \phi_R)_{n\in \mathbb{N}}$ is also bounded in the space $E_\rho$. Indeed, we have
$$\|\theta_n \phi_R\|_{E_\rho}= \sqrt{\rho}\|\theta_n \phi_R\|_{\dot{H}^1(\R)}+\|\theta_n \phi_R\|_{\dot{H}^{\frac{\alpha}{2}}(\R)}= \sqrt{\rho}\|\vn(\theta_n \phi_R)\|_{L^2(\R)}+\|(-\Delta)^{\frac{\alpha}{2}}(\theta_n \phi_R)\|_{L^2(\R)},$$
and by usual Leibniz rule and the H\"older inequalities with $\frac{1}{2}=\frac{n-\alpha}{2n}+\frac{\alpha}{2n}$, we have
$$\|\theta_n \phi_R\|_{E_\rho}\leq  \sqrt{\rho}\left(\|\phi_R\|_{L^\infty(\R)}\|\vn\theta_n\|_{L^2(\R)}+\|\vn \phi_R\|_{L^\frac{2n}{\alpha}(\R)}\|\theta_n\|_{L^{\frac{2n}{n-\alpha}}(\R)}\right)+\|(-\Delta)^{\frac{\alpha}{2}}(\theta_n \phi_R)\|_{L^2(\R)},$$
and by the Sobolev embedding $\dot{H}^{\frac{\alpha}{2}}(\R)\subset L^{\frac{2n}{n-\alpha}}(\R)$ we can write
$$\|\theta_n \phi_R\|_{E_\rho}\leq  \sqrt{\rho}\left(C_R\|\theta_n\|_{\dot{H}^1(\R)}+C_R\|\theta_n\|_{\dot{H}^{\frac{\alpha}{2}}(\R)}\right)+\|(-\Delta)^{\frac{\alpha}{2}}(\theta_n \phi_R)\|_{L^2(\R)},$$
and since $\|\cdot\|_{\dot{H}^{\frac{\alpha}{2}}(\R)}\leq \|\cdot\|_{E_\rho}$ and $\|\cdot\|_{\dot{H}^{1}(\R)}\leq \frac{1}{\sqrt{\rho}}\|\cdot\|_{E_\rho}$ we have
$$\|\theta_n \phi_R\|_{E_\rho}\leq  C_R\left(\|\theta_n\|_{E_\rho}+\sqrt{\rho}\|\theta_n\|_{E_\rho}\right)+\|(-\Delta)^{\frac{\alpha}{2}}(\theta_n \phi_R)\|_{L^2(\R)}.$$
For the second term of the right-hand side above we will use the fractional Leibniz rule and we obtain
\begin{eqnarray*}
\|\theta_n \phi_R\|_{E_\rho}&\leq& C_R\left(1+\sqrt{\rho}\right)\|\theta_n\|_{E_\rho}+\left(\|(-\Delta)^{\frac{\alpha}{2}}\theta_n\|_{L^2(\R)}\|\phi_R\|_{L^\infty(\R)}+\|\theta_n\|_{L^{\frac{2n}{n-\alpha}}(\R)}\|(-\Delta)^{\frac{\alpha}{2}}\phi_R\|_{L^{\frac{2n}{\alpha}}(\R)}\right)\\
&\leq &C_R\left(1+\sqrt{\rho}\right)\|\theta_n\|_{E_\rho}+C_R\left(\|(-\Delta)^{\frac{\alpha}{2}}\theta_n\|_{L^2(\R)}+\|\theta_n\|_{L^{\frac{2n}{n-\alpha}}(\R)}\right),
\end{eqnarray*}
using again the embedding $\dot{H}^{\frac{\alpha}{2}}(\R)\subset L^{\frac{2n}{n-\alpha}}(\R)$ and the definition of the Sobolev space $\dot{H}^{\frac{\alpha}{2}}(\R)$, we obtain
\begin{eqnarray*}
\|\theta_n \phi_R\|_{E_\rho}&\leq &C_R\left(1+\sqrt{\rho}\right)\|\theta_n\|_{E_\rho}+C_R\|\theta_n\|_{\dot{H}^{\frac{\alpha}{2}}(\R)}\\
&\leq &C_R\left(1+\sqrt{\rho}\right)\|\theta_n\|_{E_\rho}+C_R\|\theta_n\|_{E_\rho}<+\infty,
\end{eqnarray*}
and we conclude that the sequence $(\theta_n \phi_R)_{n\in \mathbb{N}}$ is also bounded in $E_\rho$.\\

Now, with this remark at hand and for a fixed $R>0$, we can assume without loss of generality that, for all $n\in \mathbb{N}$, we have $supp(\theta_n\phi_R)\subset B(0,4R)$. Since we are working in the space $\dot{H}^{\frac{\alpha}{2}}(\R)$ with $0<\alpha<2$, by the Rellich-Kondrashov lemma (see Theorem 9.16 of \cite{Brezis}), there exists a subsequence $(\theta_{n_k} \phi_R )_{k\in \mathbb{N}}$ that converges strongly in $L^p_{loc}(\R)$ for $1\le p <\frac{2n}{n-\alpha}$. This argument will lead us to the wished compactness of the application $\mathbb{T}_{R, \rho, \epsilon}$ as long as we can estimate $\|\mathbb{T}_{R, \rho, \epsilon}(\theta_n)\|_{E_\rho}$ in terms of suitable $L^p_{loc}$ norms.\\

Then, from (\ref{Continuite0}) we can write
\begin{eqnarray}
\|\mathbb{T}_{R, \rho, \epsilon}(\theta_n)\|_{E_{\rho}}&=&\sqrt{\rho}\left\|\mathfrak{L}_{\rho}\left((\mathbb{A}_{[(\theta_n \phi_R)\ast \varphi_\epsilon]}\ast \varphi_\epsilon)\phi_R\cdot \vn (\theta_n \phi_R) +f\right)\right\|_{\dot{H}^1(\R)}\notag\\
&&+\left\|\mathfrak{L}_{\rho}\left((\mathbb{A}_{[(\theta_n \phi_R)\ast \varphi_\epsilon]}\ast \varphi_\epsilon)\phi_R\cdot \vn (\theta_n \phi_R) +f\right)\right\|_{\dot{H}^{\frac{\alpha}{2}}(\R)}\notag\\[2mm]
&\leq & \sqrt{\rho}\underbrace{\left\|\mathfrak{L}_{\rho}\left((\mathbb{A}_{[(\theta_n \phi_R)\ast \varphi_\epsilon]}\ast \varphi_\epsilon)\phi_R\cdot \vn (\theta_n \phi_R)\right)\right\|_{\dot{H}^1}(\R)}_{(\alpha)}+\sqrt{\rho}\left\|\mathfrak{L}_{\rho}\left(f\right)\right\|_{\dot{H}^1(\R)}\qquad \label{Compacite0}\\
&&+\underbrace{\left\|\mathfrak{L}_{\rho}\left((\mathbb{A}_{[(\theta_n \phi_R)\ast \varphi_\epsilon]}\ast \varphi_\epsilon)\phi_R\cdot \vn (\theta_n \phi_R)\right)\right\|_{\dot{H}^{\frac{\alpha}{2}}(\R)}}_{(\beta)}+\left\|\mathfrak{L}_{\rho}\left(f\right)\right\|_{\dot{H}^{\frac{\alpha}{2}}(\R)},\label{Compacite00}
\end{eqnarray}
since the terms that contain the external force do not intervene in the study of the compactness and since we have already proven that the quantities $\left\|\mathfrak{L}_{\rho}\left(f\right)\right\|_{\dot{H}^1(\R)}$ and $\left\|\mathfrak{L}_{\rho}\left(f\right)\right\|_{\dot{H}^{\frac{\alpha}{2}}(\R)}$ are bounded by the information available over $f$ (indeed, we just proved the controls $\left\|\mathfrak{L}_{\rho}\left(f\right)\right\|_{\dot{H}^1(\R)}\leq \frac{C}{\rho}\|f\|_{\dot{H}^{-\frac{\alpha}{2}}(\R)}$ and $\left\|\mathfrak{L}_{\rho}\left(f\right)\right\|_{\dot{H}^{\frac{\alpha}{2}}(\R)}\leq \frac{C}{\rho}\|f\|_{\dot{H}^{-\frac{\alpha}{2}}(\R)}$), we will focus our study in the terms $(\alpha)$ and $(\beta)$ of the expression above. 

\begin{itemize}
\item For the term $(\alpha)$ in (\ref{Compacite0}), using the divergence free condition for the quantity $\mathbb{A}_{[(\theta \phi_R)\ast \varphi_\epsilon]}\ast \varphi_\epsilon$ we remark that we have the following identity 
\begin{eqnarray}
(\mathbb{A}_{[(\theta_n \phi_R)\ast \varphi_\epsilon]}\ast \varphi_\epsilon)\phi_R\cdot \vn (\theta_n \phi_R)&=& div\left(\mathbb{A}_{[(\theta_n \phi_R)\ast \varphi_\epsilon]}\ast \varphi_\epsilon \phi_R^2\theta_n\right)\label{DecompositionCompact}\\
&&-\mathbb{A}_{[(\theta_n \phi_R)\ast \varphi_\epsilon]}\ast \varphi_\epsilon \cdot (\vn \phi_R)(\theta_n\phi_R),\notag
\end{eqnarray}
so we can write 
\begin{eqnarray*}
\left\|\mathfrak{L}_{\rho}\left((\mathbb{A}_{[(\theta_n \phi_R)\ast \varphi_\epsilon]}\ast \varphi_\epsilon)\phi_R\cdot \vn (\theta_n \phi_R)\right)\right\|_{\dot{H}^1(\R)}&\leq &\underbrace{\left\|\mathfrak{L}_{\rho}\left(div\left(\mathbb{A}_{[(\theta_n \phi_R)\ast \varphi_\epsilon]}\ast \varphi_\epsilon \phi_R^2\theta_n\right)\right)\right\|_{\dot{H}^1(\R)}}_{(\alpha_1)}\\
&&+\underbrace{\left\|\mathfrak{L}_{\rho}\left(\mathbb{A}_{[(\theta_n \phi_R)\ast \varphi_\epsilon]}\ast \varphi_\epsilon \cdot (\vn \phi_R)(\theta_n\phi_R)\right)\right\|_{\dot{H}^1(\R)}}_{(\alpha_2)}.
\end{eqnarray*}
Now for the term $(\alpha_1)$ in the right-hand side above, by the definition of the operator $\mathfrak{L}_{\rho}$ given in (\ref{DefopL}) and by the properties of the Sobolev spaces, we write
\begin{eqnarray}
\alpha_1&=&\left\| \frac{1}{[\rho(-\Delta)+(-\Delta)^{\frac{\alpha}{2}}]}div\left(\mathbb{A}_{[(\theta_n \phi_R)\ast \varphi_\epsilon]}\ast \varphi_\epsilon \phi_R^2\theta_n\right)\right\|_{\dot{H}^1(\R)}\notag\\
&\leq &C\left\| \frac{(-\Delta)}{[\rho(-\Delta)+(-\Delta)^{\frac{\alpha}{2}}]}\left(\mathbb{A}_{[(\theta_n \phi_R)\ast \varphi_\epsilon]}\ast \varphi_\epsilon \phi_R^2\theta_n\right)\right\|_{L^2(\R)}\notag\\
&\leq &\frac{C}{\rho}\left\|\mathbb{A}_{[(\theta_n \phi_R)\ast \varphi_\epsilon]}\ast \varphi_\epsilon \phi_R^2\theta_n\right\|_{L^2(\R)},\label{EstimationCompaciteAlpha1}
\end{eqnarray}
where we used the Lemma \ref{LemmeControlOpFourier} to deduce the last estimate. Since we have $2<\frac{2n}{n-\alpha}$, thus by the H\"older inequalities, the Young inequalities and by the boundedness properties of the operator $\mathbb{A}$, we can write
\begin{eqnarray}
\alpha_1&\leq& \frac{C}{\rho}\left\| \mathbb{A}_{[(\theta_n \phi_R)\ast \varphi_\epsilon]}\ast \varphi_\epsilon \right\|_{L^\infty(\R)}\left\|\phi_R^2\theta_n\right\|_{L^2(\R)}\leq \frac{C}{\rho}\left\| \mathbb{A}_{[(\theta_n \phi_R)\ast \varphi_\epsilon]}\right\|_{L^2(\R)}\| \varphi_\epsilon\|_{L^2(\R)}\left\|\phi_R^2\theta_n\right\|_{L^2(\R)}\notag\\
&\leq & \frac{C_{R,\mathbb{A}, \epsilon}}{\rho}\|(\theta_n \phi_R)\ast \varphi_\epsilon\|_{L^2(\R)}\|\theta_n\phi_R\|_{L^2(\R)}\leq \frac{C_{R,\mathbb{A}, \epsilon}}{\rho}\|\theta_n \phi_R\|_{L^2(\R)}\|\theta_n\phi_R\|_{L^2(\R)}\notag\\
&\leq & \frac{C_{R,\mathbb{A}, \epsilon}}{\rho}\|\theta_n \phi_R\|_{L^2(B(0,4R))}\|\theta_n\phi_R\|_{L^2(B(0,4R))},\label{EstimationCompaciteAlpha11}
\end{eqnarray}
where we used the support properties of the function $\phi_R$. Since the sequence $(\theta_n\phi_R)_{n\in \mathbb{N}}$ is bounded in the space $\dot{H}^{\frac{\alpha}{2}}(\R)$ and since $2<\frac{2n}{n-\alpha}$, we can extract a subsequence $(\theta_{n_k}\phi_R)_{k\in \mathbb{N}}$ that converges strongly in the space $L^2(B(0,4R))$ and this will provide us the wished compactness property associated to the term $(\alpha_1)$.\\

For the term $(\alpha_2)$ we proceed as follows: for some $0<s<1$ such that $\alpha<1+s<2$ we write
\begin{eqnarray*}
\alpha_2&=&\left\| \frac{(-\Delta)^{\frac{1+s}{2}}}{[\rho(-\Delta)+(-\Delta)^{\frac{\alpha}{2}}]} \frac{1}{(-\Delta)^{\frac{s}{2}}}\left(\mathbb{A}_{[(\theta_n \phi_R)\ast \varphi_\epsilon]}\ast \varphi_\epsilon \cdot (\vn \phi_R)(\theta_n\phi_R)\right)\right\|_{L^2(\R)}\\
&\leq & \frac{C}{\rho}\left\| \mathbb{A}_{[(\theta_n \phi_R)\ast \varphi_\epsilon]}\ast \varphi_\epsilon \cdot (\vn \phi_R)(\theta_n\phi_R)\right\|_{\dot{H}^{-s}(\R)}\\
&\leq & \frac{C}{\rho}\left\| \mathbb{A}_{[(\theta_n \phi_R)\ast \varphi_\epsilon]}\ast \varphi_\epsilon \cdot (\vn \phi_R)(\theta_n\phi_R)\right\|_{L^{\frac{2n}{n+2s}}},
\end{eqnarray*}
where we applied the Lemma \ref{LemmeControlOpFourier} and we used the definition of the Sobolev space $\dot{H}^{-s}(\R)$ as well as the embedding $L^{\frac{2n}{n+2s}}(\R)\subset \dot{H}^{-s}(\R)$. Now, by the H\"older inequalities with $\frac{n+2s}{2n}=\frac{1}{2}+\frac{s}{n}$ and with the Young inequalities with $1+\frac{s}{n}=\frac{1}{2}+\frac{n+2s}{2n}$, we have
\begin{eqnarray*}
\alpha_2&\leq &\frac{C}{\rho}\left\|  \mathbb{A}_{[(\theta \phi_R)\ast \varphi_\epsilon]}\ast \varphi_\epsilon\right\|_{L^{\frac{n}{s}}(\R)} \left\|\vn \phi_R\right\|_{L^{\infty}(\R)} \left\|\theta_n\phi_R\right\|_{L^{2}(\R)}\\
&\leq & \frac{C_{R}}{\rho}\|\mathbb{A}_{[(\theta \phi_R)\ast \varphi_\epsilon]}\|_{L^2(\R)}\|\varphi_\epsilon\|_{L^{\frac{2n}{n+2s}}(\R)} \|\theta_n\phi_R\|_{L^2(\R)}\\
&\leq & \frac{C_{R, \mathbb{A}, \varepsilon}}{\rho}\|(\theta \phi_R)\ast \varphi_\epsilon\|_{L^2(\R)}\|\theta_n\phi_R\|_{L^2(\R)}\leq  \frac{C_{R, \mathbb{A}, \varepsilon}}{\rho}\|\theta \phi_R\|_{L^2(\R)}\|\theta_n\phi_R\|_{L^2(\R)}\\
&\leq& \frac{C_{R, \mathbb{A}, \varepsilon}}{\rho}\|\theta \phi_R\|_{L^2(B(0,4R))}\|\theta_n\phi_R\|_{L^2(B(0,4R))},
\end{eqnarray*}
where we used the boundedness of the operator $\mathbb{A}$ and the support properties of the function $\phi_R$. With this control we can deduce as before the wished compactness property associated to the term $(\alpha_2)$.\\

With these compactness properties for the terms $(\alpha_1)$ and $(\alpha_2)$, we can easily deduce the compactness of the quantity $(\alpha)$ given in (\ref{Compacite0}).\\[3mm]

\item For the term $(\beta)$ in (\ref{Compacite00}) we have, using the decomposition (\ref{DecompositionCompact}):
\begin{eqnarray*}
\left\|\mathfrak{L}_{\rho}\left((\mathbb{A}_{[(\theta_n \phi_R)\ast \varphi_\epsilon]}\ast \varphi_\epsilon)\phi_R\cdot \vn (\theta_n \phi_R)\right)\right\|_{\dot{H}^{\frac{\alpha}{2}}(\R)}\leq \underbrace{\left\|\mathfrak{L}_{\rho}\left(div\left(\mathbb{A}_{[(\theta_n \phi_R)\ast \varphi_\epsilon]}\ast \varphi_\epsilon \phi_R^2\theta_n\right)\right)\right\|_{\dot{H}^{\frac{\alpha}{2}}(\R)}}_{(\beta_1)}\\
+\underbrace{\left\|\mathfrak{L}_{\rho}\left(\mathbb{A}_{[(\theta_n \phi_R)\ast \varphi_\epsilon]}\ast \varphi_\epsilon \cdot (\vn \phi_R)(\theta_n\phi_R)\right)\right\|_{\dot{H}^{\frac{\alpha}{2}}(\R)}}_{(\beta_2)},
\end{eqnarray*}
and we will study these terms separately.\\ 

For the quantity $(\beta_1)$ above we write, by the definition of the operator $\mathfrak{L}_{\rho}$ given in (\ref{DefopL}) and by the properties of the Sobolev spaces:
\begin{eqnarray*}
\beta_1&=&\left\| \frac{1}{[\rho(-\Delta)+(-\Delta)^{\frac{\alpha}{2}}]} div\left(\mathbb{A}_{[(\theta_n \phi_R)\ast \varphi_\epsilon]}\ast \varphi_\epsilon \phi_R^2\theta_n\right)\right\|_{\dot{H}^{\frac{\alpha}{2}}(\R)}\\
&\leq &\left\| \frac{(-\Delta)^{\frac{1+\frac{\alpha}{2}}{2}}}{[\rho(-\Delta)+(-\Delta)^{\frac{\alpha}{2}}]}\left(\mathbb{A}_{[(\theta_n \phi_R)\ast \varphi_\epsilon]}\ast \varphi_\epsilon \phi_R^2\theta_n\right)\right\|_{L^2(\R)},
\end{eqnarray*}
since $\frac{\alpha}{2}<\frac{1+\frac{\alpha}{2}}{2}<1$, by the Lemma \ref{LemmeControlOpFourier}, we have the estimate
$$\beta_1\leq \frac{C}{\rho}\left\|\mathbb{A}_{[(\theta_n \phi_R)\ast \varphi_\epsilon]}\ast \varphi_\epsilon \phi_R^2\theta_n\right\|_{L^2(\R)},$$
at this point, following the same arguments as above (see (\ref{EstimationCompaciteAlpha1})-(\ref{EstimationCompaciteAlpha11})), we obtain the control
$$\beta_1\leq  \frac{C_{R,\mathbb{A}, \epsilon}}{\rho}\|\theta_n \phi_R\|_{L^2(B(0,4R))}\|\theta_n\phi_R\|_{L^2(B(0,4R))},$$
from which we deduce the compactness property associated to the term $(\beta_1)$.\\

For the term $(\beta_2)$, we have as before:
\begin{eqnarray*}
\beta_2&=&\left\| \frac{(-\Delta)^{\frac{\alpha}{2}}}{[\rho(-\Delta)+(-\Delta)^{\frac{\alpha}{2}}]} \frac{1}{(-\Delta)^{\frac{\mathfrak{\alpha}}{4}}}\left(\mathbb{A}_{[(\theta_n \phi_R)\ast \varphi_\epsilon]}\ast \varphi_\epsilon \cdot (\vn \phi_R)(\theta_n\phi_R)\right)\right\|_{L^2(\R)}\\
&\leq& \frac{C}{\rho}\left\|\mathbb{A}_{[(\theta_n \phi_R)\ast \varphi_\epsilon]}\ast \varphi_\epsilon \cdot (\vn \phi_R)(\theta_n\phi_R)\right\|_{\dot{H}^{-{\frac{\alpha}{2}}}(\R)}\\
&\leq &\frac{C}{\rho}\left\|\mathbb{A}_{[(\theta_n \phi_R)\ast \varphi_\epsilon]}\ast \varphi_\epsilon \cdot (\vn \phi_R)(\theta_n\phi_R)\right\|_{L^{\frac{2n}{n+\alpha}}(\R)},
\end{eqnarray*}
where we used the Lemma \ref{LemmeControlOpFourier} and the Hardy-Littlewood-Sobolev inequality.
Now, by the H\"older inequality with $\frac{n+\alpha}{2n}=\frac{1}{2}+\frac{\alpha}{2n}$, we can write
\begin{eqnarray*}
\beta_2&\leq &\frac{C}{\rho}\left\|\mathbb{A}_{[(\theta_n \phi_R)\ast \varphi_\epsilon]}\ast \varphi_\epsilon\right\|_{L^{\frac{2n}{\alpha}}(\R)}\left\| \vn \phi_R\right\|_{L^{\infty}(\R)}\left\|\theta_n\phi_R\right\|_{L^{2}(\R)}\\
&\leq & \frac{C}{\rho}\left\|\mathbb{A}_{[(\theta_n \phi_R)\ast \varphi_\epsilon]}\right\|_{L^{2}(\R)}\| \varphi_\epsilon\|_{L^{\frac{2n}{n+\alpha}}(\R)}\left\| \vn \phi_R\right\|_{L^{\infty}(\R)}\left\|\theta_n\phi_R\right\|_{L^{2}(\R)},
\end{eqnarray*}
where we used the Young inequalities with $1+\frac{\alpha}{2n}=\frac{1}{2}+\frac{n+\alpha}{2n}$. Since the operator $\mathbb{A}$ is bounded, and by the properties of the function $\phi_R$, we have
\begin{eqnarray*}
\beta_2&\leq &\frac{C_{R, \epsilon}}{\rho}\left\|\mathbb{A}_{[(\theta_n \phi_R)\ast \varphi_\epsilon]}\right\|_{L^{2}(\R)}\left\|\theta_n\phi_R\right\|_{L^{2}(\R)}\leq \frac{C_{R, \mathbb{A}, \epsilon}}{\rho}\left\|(\theta_n \phi_R)\ast \varphi_\epsilon\right\|_{L^{2}(\R)}\left\|\theta_n\phi_R\right\|_{L^{2}(\R)}\\
&\leq &\frac{C_{R, \mathbb{A}, \epsilon}}{\rho}\|\theta_n\phi_R\|_{L^2(B(0,4R))}\|\theta_n\phi_R\|_{L^2(B(0,4R))},
\end{eqnarray*}
from which we deduce the compactness property associated to the term $(\beta_2)$.\\

We have obtained the compactness property for the terms $(\beta_1)$ and $(\beta_2)$ defined above and we deduce then the compactness property for the term $(\beta)$ given in (\ref{Compacite00}).\\

\end{itemize}
From this separated study of the quantities $(\alpha)$ and $(\beta)$, we finally obtain the compactness of the operator $\mathbb{T}_{R, \rho, \epsilon}$ in the space $E_\rho$ and this ends the proof of the Proposition \ref{PropoCompacidad}. \hfill $\blacksquare$\\

\begin{Proposition}[A priori estimates]\label{PropoApriori}
In the setting of the space $E_\rho$, if $\theta=\lambda \mathbb{T}_{R, \rho, \epsilon}(\theta)$ for any $\lambda\in [0,1]$, then we have the control $\|\theta\|_{E_{\rho}}\leq M$.
\end{Proposition}
{\bf Proof.} If we have for all $0\leq \lambda\leq 1$ the identity
$$\theta=\lambda \mathbb{T}_{R, \rho, \epsilon}(\theta),$$
by the definition of the operator $\mathbb{T}_{R, \rho, \epsilon}$ we can write
$$\theta=\lambda\left[\frac{1}{[\rho(-\Delta)+(-\Delta)^{\frac{\alpha}{2}}]}\left((\mathbb{A}_{[(\theta \phi_R)\ast \varphi_\epsilon]}\ast \varphi_\epsilon)\phi_R\cdot \vn (\theta \phi_R) +f\right)\right],$$
we then obtain the expression
$$[\rho(-\Delta)+(-\Delta)^{\frac{\alpha}{2}}]\theta=\lambda\left((\mathbb{A}_{[(\theta \phi_R)\ast \varphi_\epsilon]}\ast \varphi_\epsilon)\phi_R\cdot \vn (\theta \phi_R)\right) +\lambda f,$$
and we have
$$\rho\int_{\R}[(-\Delta)\theta]\theta dx+\int_{\R}[(-\Delta)^{\frac{\alpha}{2}}\theta] \theta dx=\lambda\int_{\R}\left((\mathbb{A}_{[(\theta \phi_R)\ast \varphi_\epsilon]}\ast \varphi_\epsilon)\phi_R\cdot \vn (\theta \phi_R)\right)\theta dx +\lambda \int_{\R}f\theta dx.$$
Since $div(\mathbb{A}_{[(\theta \phi_R)\ast \varphi_\epsilon]}\ast \varphi_\epsilon)=0$, we obtain by an integration by parts that 
$$\int_{\R}\left((\mathbb{A}_{[(\theta \phi_R)\ast \varphi_\epsilon]}\ast \varphi_\epsilon)\phi_R\cdot \vn (\theta \phi_R)\right)\theta dx=-\int_{\R}\left((\mathbb{A}_{[(\theta \phi_R)\ast \varphi_\epsilon]}\ast \varphi_\epsilon)\phi_R\cdot \vn (\theta \phi_R)\right)\theta dx,$$
from which we deduce that 
$$\int_{\R}\left((\mathbb{A}_{[(\theta \phi_R)\ast \varphi_\epsilon]}\ast \varphi_\epsilon)\phi_R\cdot \vn (\theta \phi_R)\right)\theta dx=0,$$
so we obtain the identity 
$$\rho\|\theta\|_{\dot{H}^1(\R)}^2+\|\theta\|_{\dot{H}^{\frac{\alpha}{2}}(\R)}^2=\lambda \int_{\R}f\theta dx.$$
Now, by the $\dot{H}^{\frac{\alpha}{2}}-\dot{H}^{-\frac{\alpha}{2}}$ duality we can write
$$\rho\|\theta\|_{\dot{H}^1(\R)}^2+\|\theta\|_{\dot{H}^{\frac{\alpha}{2}}(\R)}^2\leq \lambda \|f\|_{\dot{H}^{-\frac{\alpha}{2}}(\R)}\|\theta\|_{\dot{H}^{\frac{\alpha}{2}}(\R)},$$
and by the Young inequalities we have
$$\rho\|\theta\|_{\dot{H}^1(\R)}^2+\|\theta\|_{\dot{H}^{\frac{\alpha}{2}}(\R)}^2\leq \frac{\lambda}{2} \|f\|_{\dot{H}^{-\frac{\alpha}{2}}(\R)}^2+\frac{\lambda}{2}\|\theta\|_{\dot{H}^{\frac{\alpha}{2}}(\R)}^2,$$
which can be rewritten as
$$\rho\|\theta\|_{\dot{H}^1(\R)}^2+\left(1-\frac{\lambda}{2}\right)\|\theta\|_{\dot{H}^{\frac{\alpha}{2}}(\R)}^2\leq \frac{\lambda}{2} \|f\|_{\dot{H}^{-\frac{\alpha}{2}}(\R)}^2.$$
Finally, since $0\leq \lambda \leq 1$ and since $f\in \dot{H}^{-\frac{\alpha}{2}}(\R)$ by hypothesis, we can write 
\begin{equation}\label{EstimationUniformes}
\rho\|\theta\|_{\dot{H}^1(\R)}^2+\|\theta\|_{\dot{H}^{\frac{\alpha}{2}}(\R)}^2\leq \|f\|_{\dot{H}^{-\frac{\alpha}{2}}(\R)}^2,\end{equation}
since $f\in \dot{H}^{-\frac{\alpha}{2}}(\R)$ by hypothesis, we deduce the wished estimate, 
$$\|\theta\|_{E_\rho}\leq M<+\infty,$$
and the proof of the Proposition \ref{PropoApriori} is complete. \hfill $\blacksquare$
\begin{Remarque}
Note that the estimate (\ref{EstimationUniformes}) gives a uniform control (in $0<\rho, \epsilon<1$ and in $R>1$) of the quantity $\|\theta\|_{\dot{H}^{\frac{\alpha}{2}}(\R)}$.
\end{Remarque}
\subsection*{End of the proof of the Theorem \ref{Theo_ExistenceStationnaires}.}
Gathering Propositions \ref{PropoContinuidad}, \ref{PropoCompacidad} and \ref{PropoApriori}, we can apply the Theorem \ref{Teo_Schaefer} to obtain a solution $\theta=\theta_{\rho, \epsilon, R}$ of the equation (\ref{EquationRegulariseeStationnaire}). Now, in order to recover a solution of the original equation (\ref{SQGStationnaire}) we need to take the limits $\rho, \epsilon\to 0$ and $R\to +\infty$. To do so, we will exploit the uniform estimate (\ref{EstimationUniformes}), indeed, 
since $\theta$ is uniformly bounded by (\ref{EstimationUniformes}) in the space $\dot{H}^{\frac{\alpha}{2}}(\R)$, by the Banach-Alaoglu theorem there exists a limit $\theta_{\rho, \epsilon, R}\underset{\rho, \epsilon\to 0}{\longrightarrow}\theta_R$ in the weak-$*$ sense. Furthermore we have $(-\Delta)^{\frac{\alpha}{2}}\theta_{\rho, \epsilon, R}\underset{\rho, \epsilon\to 0}{\longrightarrow}(-\Delta)^{\frac{\alpha}{2}}\theta_R$ in the sense of tempered distributions. For the non linear term $(\mathbb{A}_{[(\theta \phi_R)\ast \varphi_\epsilon]}\ast \varphi_\epsilon)\phi_R\cdot \vn (\theta \phi_R)$, we will need some local strong convergence in a suitable Lebesgue space, which can be obtained by the Rellich-Kondrashov theorem since we have the uniform control (\ref{EstimationUniformes}). We can thus obtain, in a weak sense the limit $(\mathbb{A}_{[(\theta_{\rho, \epsilon, R} \phi_R)\ast \varphi_\epsilon]}\ast \varphi_\epsilon)\phi_R\cdot \vn (\theta_{\rho, \epsilon, R} \phi_R)\underset{\rho, \epsilon\to 0}{\longrightarrow}(\mathbb{A}_{[\theta_{R} \phi_R]}\phi_R)\cdot \vn (\theta_{R} \phi_R)$. These arguments are classical for studying nonlinear terms. Now we must consider the limit when $R\to+\infty$ and for this we can proceed in the same manner as before since the control (\ref{EstimationUniformes}) is also uniform in $R$. Thus, with this arguments, we obtain a solution $\theta\in \dot{H}^{\frac{\alpha}{2}}(\R)$ which is a weak solution of the equation
$$(-\Delta)^{\frac \alpha 2}\theta-\mathbb{A}_{[\theta]}\cdot \vn\theta-f=0,$$
and now the proof of the Theorem \ref{Theo_ExistenceStationnaires} is complete.\hfill $\blacksquare$
\mysection{Regularity}\label{Secc_Regularity}
The problem of the regularity of weak solutions to the system (\ref{SQGStationnaire}) is of course related to the power of the fractional Laplacian $0<\alpha<2$. Our starting point is a function $\theta \in \dot{H}^{\frac{\alpha}{2}}(\R)$ that satisfies in the weak sense the equation (\ref{SQGStationnaire}) where, for the sake of simplicity we will assume from now on that the external force $f$ is null.\\

We start with the case $\frac{n+2}{3}<\alpha<2$ and we rewrite the equation (\ref{SQGStationnaire}) as follows:
$$\theta=\frac{1}{(-\Delta)^{\frac \alpha 2}}\mathbb{A}_{[\theta]}\cdot \vn\theta=\frac{1}{(-\Delta)^{\frac \alpha 2}}div(\mathbb{A}_{[\theta]}\theta),$$
where we used the divergence free property of $\mathbb{A}$. Now, for some index $\sigma>0$ that will be defined later, we write
$$\|(-\Delta)^{\frac{\sigma}{2}}\theta\|_{L^2(\R)}=\left\|(-\Delta)^{\frac{\sigma-\alpha}{2}}div(\mathbb{A}_{[\theta]}\theta)\right\|_{L^2(\R)}\leq C\|\mathbb{A}_{[\theta]}\theta\|_{\dot{H}^{\sigma-\alpha+1}(\R)}.$$
At this point we apply the product law given in the Lemma \ref{ProductRule} to obtain (since $\theta\in \dot{H}^{\frac{\alpha}{2}}(\R)$ and since the operator $\mathbb{A}$ is bounded in Sobolev spaces):
\begin{equation}\label{Application_LoiProduit}
\|\mathbb{A}_{[\theta]}\theta\|_{\dot{H}^{\sigma-\alpha+1}(\R)}\leq C\|\mathbb{A}_{[\theta]}\|_{\dot{H}^{\frac{\alpha}{2}}(\R)}\|\theta\|_{\dot{H}^{\frac{\alpha}{2}}(\R)}\leq C\|\theta\|_{\dot{H}^{\frac{\alpha}{2}}(\R)}\|\theta\|_{\dot{H}^{\frac{\alpha}{2}}(\R)}<+\infty,
\end{equation}
as long as $\sigma-\alpha+1=\alpha-\frac{n}{2}$, from which we deduce that $\sigma=2\alpha-\frac{n}{2}-1$. This is a gain of regularity with respect to the information $\theta\in \dot{H}^{\frac{\alpha}{2}}(\R)$ as long as we have $\sigma=2\alpha-\frac{n}{2}-1>\frac{\alpha}{2}$, \emph{i.e.} if we have $\alpha>\frac{n+2}{3}$ which is the case studied here. Once we obtain this first gain of regularity, by iterating this process we easily obtain that the solutions of the equation (\ref{SQGStationnaire}) are smooth. \\

We consider now the case $1<\alpha\leq \frac{n+2}{3}$. Here the general framework $\theta\in \dot{H}^{\frac{\alpha}{2}}(\R)$ seems to be not enough to obtain a gain of regularity when applying the Lemma \ref{ProductRule} in the estimate (\ref{Application_LoiProduit}) above. For simplicity we will use here an additional hypothesis given by $\theta\in L^\infty(\R)$ and instead of Lemma \ref{ProductRule} we use the Leibniz fractional inequality given in Lemma \ref{FracLeibniz}. Thus, instead of (\ref{Application_LoiProduit}) we write:
$$\|\mathbb{A}_{[\theta]}\theta\|_{\dot{H}^{\sigma-\alpha+1}(\R)}\leq C(\|\mathbb{A}_{[\theta]}\|_{\dot{H}^{\sigma-\alpha+1}(\R)}\|\theta\|_{L^\infty(\R)}+\|\mathbb{A}_{[\theta]}\|_{L^\infty(\R)}\|\theta\|_{\dot{H}^{\sigma-\alpha+1}(\R)}).$$
We will assume moreover that the operator $\mathbb{A}$ is also bounded in the space $L^\infty(\R)$, so we can write
$$\|\mathbb{A}_{[\theta]}\theta\|_{\dot{H}^{\sigma-\alpha+1}(\R)}\leq C\|\theta\|_{\dot{H}^{\sigma-\alpha+1}(\R)}\|\theta\|_{L^\infty(\R)},$$
which is a finite quantity as long as $\sigma-\alpha+1=\frac{\alpha}{2}$, which gives $\sigma=\frac{3}{2}\alpha-1$. But since $1<\alpha$ we have $\sigma>\frac{\alpha}{2}$ and we have obtained a gain of regularity as we have proved that $\theta\in \dot{H}^{\sigma}(\R)$. Again, by iteration we obtain that the solutions of the equation (\ref{SQGStationnaire}) are smooth.\\

The case $0<\alpha\leq 1$ is much harder to study and, to the best of our knowledge, it is still an open problem. \\


\noindent {\bf Acknowledgment.} The work of the second author has been partially supported by the program MATH-AMSUD 23-MATH-18.


\end{document}